\documentclass[10pt]{amsart}

\usepackage{amsmath,amsthm,amssymb}
\usepackage{epsfig}
\usepackage[margin=1.3in]{geometry}
\usepackage{lscape}

\begin{document}

%%%%%%%%%%%%%%%%%%%%%%%%%%%%%%%%%%%%%%%%%%%%%%%%%%%%%%%%%%%%%%%%%%%%%%%%%%%%%%%%%%%%%
%  Group Notation
%%%%%%%%%%%%%%%%%%%%%%%%%%%%%%%%%%%%%%%%%%%%%%%%%%%%%%%%%%%%%%%%%%%%%%%%%%%%%%%%%%%%%

\newcommand{\tot}  {{\Omega}}                           %operad

\newcommand{\Op}  {{\mathcal O}}                           %operad
\newcommand{\Opu}  {{\mathcal U}}                           %operad
\newcommand{\Opv}  {{\mathcal V}}                           %operad
\newcommand{\PGL} {\Pj\Gl_2(\R)}                           %PGL_2(R)
\newcommand{\PGLC} {\Pj\Gl_2(\Cx)}                         %PGL_2(C)
\newcommand{\RP} {\R\Pj^1}                                 %RP^1
\newcommand{\CP} {\Cx\Pj^1}                                %CP^1
\newcommand{\Cx} {{\mathbb C}}                             %complex
\newcommand{\R} {{\mathbb R}}                              %reals
\newcommand{\I} {{\mathbb I}}                              %interval
\newcommand{\Z} {{\mathbb Z}}                              %integers
\newcommand{\Pj} {{\mathbb P}}                             %projective space
\newcommand{\T} {{\mathbb T}}                              %torus
\newcommand{\Sg} {{\mathbb S}}                             %circle
\newcommand{\Gl} {{\rm Gl}}                                %g.linear group

\newcommand{\sg} {\sigma}
\newcommand{\sghat} {\hat{\sg}}

\newcommand{\suchthat} {\:\: | \:\:}
\newcommand{\ore} {\ \ {\it or} \ \ }
\newcommand{\oand} {\ \ {\it and} \ \ }

\newcommand{\At} {W_{\mathcal A}}
\newcommand{\Sp} {W_{\mathcal S}}
\newcommand{\Eu} {W_{\mathcal E}}
\newcommand{\Qs} {\mathcal {Q_S}}
\newcommand{\Qe} {\mathcal {Q_E}}

%%%%%%%%%%%%%%%%%%%%%%%%%%%
% Configuration spaces
%%%%%%%%%%%%%%%%%%%%%%%%%%%

\newcommand{\Con} {\mathbb F}

\newcommand{\Bc} [2] {{B_{#2}\langle{#1}\rangle}}                         %closed  unordered
\newcommand{\B}  [2] {{B_{#2}({#1})}}                                     %open    unordered
\newcommand{\BC} [2] {{B_{#2}[{#1}]}}                                     %compact unordered
\newcommand{\F}  [2] {{C_{#2}({#1})}}                                     %open    ordered
\newcommand{\FC} [2] {{C_{#2}\langle{#1}\rangle}}                         %closed  ordered

\newcommand{\cM} [1] {\ensuremath{{\mathcal M}_{0, #1}}}                  %open complex
\newcommand{\CM} [1] {\ensuremath{{\overline{\mathcal M}}{_{0, #1}}}}     %cpt. complex
\newcommand{\oM} [1] {\ensuremath{{\mathcal M}_{0, #1}(\R)}}              %open
\newcommand{\M} [1] {\ensuremath{{\overline{\mathcal M}}{_{0, #1}(\R)}}}  %compact DMK
\newcommand{\Aff} {{\rm Aff}}                                             %affine automorphisms

%%%%%%%%%%%%%%%%%%%%%%%%%%%
% Some Algebraic Geometry
%%%%%%%%%%%%%%%%%%%%%%%%%%%

\newcommand{\al}{\alpha}
\newcommand{\be}{\beta}
\newcommand{\ga}{\gamma}

%%%%%%%%%%%%%%%%%%%%%%%%%%%
% Hyperplane commands
%%%%%%%%%%%%%%%%%%%%%%%%%%%

\newcommand{\Hy} {{\mathcal H}}                                            %hyperplanes that fix e_k
\newcommand{\Hs} [1] {{\mathcal H}_s{#1}}                                  %hyperplanes stabilize

%%%%%%%%%%%%%%%%%%%%%%%%%%%
% Coxeter complexes and hyperplanes
%%%%%%%%%%%%%%%%%%%%%%%%%%%

\newcommand{\C} [1]  {\mathcal C{#1}}                                      %cox complex
\newcommand{\Cp} [1] {{\Pj\C{}}({#1})}                                     %Projective cox complex
\newcommand{\Cm} [1] {\C{} ({#1})_{\#}}                                    %min blowup of cox complex
\newcommand{\Cpm} [1] {\Pj{}\Cm{#1}}                                       %min blowup of Projective cox complex
\newcommand{\Min} [1] {\textrm{Min}(\C{#1})}                               %minimal building set

\newcommand{\An} [1] {A_{#1}}
\newcommand{\Bn} [1] {B_{#1}}
\newcommand{\Dn} [1] {D_{#1}}

\newcommand{\Ant} [1] {\widetilde{A}_{#1}}
\newcommand{\Bnt} [1] {\widetilde{B}_{#1}}
\newcommand{\Cnt} [1] {\widetilde{C}_{#1}}
\newcommand{\Dnt} [1] {\widetilde{D}_{#1}}

\newcommand{\cas} [1] {{\C{A_{#1}}}}                         %coxeter A_n complex
\newcommand{\cbs} [1] {{\C{B_{#1}}}}                         %coxeter B_n complex
\newcommand{\cds} [1] {{\C{D_{#1}}}}                         %coxeter D_n complex

\newcommand{\cat} [1] {{\T\C{\widetilde A_{#1}}}}            %coxeter affine A_n complex
\newcommand{\cbt} [1] {{\T\C{\widetilde B_{#1}}}}            %coxeter affine B_n complex
\newcommand{\cct} [1] {{\T\C{\widetilde C_{#1}}}}            %coxeter affine C_n complex
\newcommand{\cdt} [1] {{\T\C{\widetilde D_{#1}}}}            %coxeter affine D_n complex

\newcommand{\csds} [2] {{\C{D_{#1, #2}}}}                    %D_(n,k) complex
\newcommand{\csbt} [2] {{\T\C{\widetilde B_{#1, #2}}}}       %affine B_(n,k) complex
\newcommand{\csdt} [2] {{\T\C{\widetilde D_{#1, #2}}}}       %affine D_(n,k) complex

%%%%%%%%%%%%%%%%%%%%%%%%%%%%%%%%%%%%%%%%%%%%%%%%%%%%%%%%%%%%%%%%%%%%%%%%%%%%%%%%%%%%%

\newcommand{\Cas} [1] {{\C{(A_{#1})}_\#}}                         %cpt coxeter A_n complex
\newcommand{\Cbs} [1] {{\C{(B_{#1})}_\#}}                         %cpt coxeter B_n complex
\newcommand{\Cds} [1] {{\C{(D_{#1})}_\#}}                         %cpt coxeter D_n complex

\newcommand{\Cat} [1] {{\T\C{(\widetilde A_{#1})}_\#}}            %cpt affine A_n complex
\newcommand{\Cbt} [1] {{\T\C{(\widetilde B_{#1})}_\#}}            %cpt affine B_n complex
\newcommand{\Cct} [1] {{\T\C{(\widetilde C_{#1})}_\#}}            %cpt affine C_n complex
\newcommand{\Cdt} [1] {{\T\C{(\widetilde D_{#1})}_\#}}            %cpt affine D_n complex

\newcommand{\Csds} [2] {{\C{(D_{#1, #2})}_\#}}                    %cpt D_(n,k) complex
\newcommand{\Csbt} [2] {{\T\C{(\widetilde B_{#1, #2})}_\#}}       %cpt affine B_(n,k) cpx
\newcommand{\Csdt} [2] {{\T\C{(\widetilde D_{#1, #2})}_\#}}       %cpt affine D_(n,k) cpx

%%%%%%%%%%%%%%%%%%%%%%%%%%%
% Graph commands
%%%%%%%%%%%%%%%%%%%%%%%%%%%

\newcommand{\Pol}{{\mathcal{P}}}                          %Polytope
\newcommand{\Cox}{{\Gamma}}                               %coxeter graph
\newcommand{\Coxrec}[1]{{\Gamma^*_{#1}}}                  %reconnected complement
\newcommand{\f}[1]{{X^{a}_{#1}}}                          %fish tail
\newcommand{\ff}[1]{{X^{b}_{#1}}}                         %forked fish tail
\newcommand{\df}[1]{{X^{c}_{#1}}}                         %double fish tail
\newcommand{\Fp}{X_4}                                     %4-permutahedron

%%%%%%%%%%%%%%%%%%%%%%%%%%%
% Polyhedra commands
%%%%%%%%%%%%%%%%%%%%%%%%%%%

\newcommand{\PFp}{{\mathcal{P} \Fp}}                  %polytope of type F_4
\newcommand{\Pf}[1]{{\mathcal{P}X^{a}_{#1}}}          %polytope of fish tail
\newcommand{\Pft}[1]{{\mathcal{P}X^{b}_{#1}}}         %polytope of forked fish tail
\newcommand{\Pdf}[1]{{\mathcal{P}X^{c}_{#1}}}         %polytope of double fish tail
\newcommand{\PA}[1]{{\mathcal{P}A_{#1}}}              %polytope of A
\newcommand{\PB}[1]{{\mathcal{P}B_{#1}}}              %polytope of B
\newcommand{\PD}[1]{{\mathcal{P}D_{#1}}}              %polytope of D
\newcommand{\PAt}[1]{{\mathcal{P}\widetilde A_{#1}}}  %polytope of Atilde
\newcommand{\PBt}[1]{{\mathcal{P}\widetilde B_{#1}}}  %polytope of Btilde
\newcommand{\PCt}[1]{{\mathcal{P}\widetilde C_{#1}}}  %polytope of Ctilde
\newcommand{\PDt}[1]{{\mathcal{P}\widetilde D_{#1}}}  %polytope of Dtilde
\newcommand{\PG}[1]{{\mathcal{P}\Gamma_{#1}}}         %polytope of gamma
\newcommand{\PGre}[1]{{\mathcal{P}\Gamma^*_{#1}}}     %polytope of gammareconnected

%%%%%%%%%%%%%%%%%%%%%%%%%%%
% Config Spaces commands
%%%%%%%%%%%%%%%%%%%%%%%%%%%
\newcommand{\Acon}[1]{\F{\R}{#1}}
\newcommand{\Bcon}[1]{\F{\R_{\bullet}}{\bar{#1}}}
\newcommand{\Dcon}[1]{\F{\R_{\circ}}{\bar{#1}}}
\newcommand{\Atcon}[1]{\F{\Sg^\prime}{#1}}
\newcommand{\Btcon}[1]{\F{\Sg^{\bullet}_{\circ}}{\bar{#1}}}
\newcommand{\Ctcon}[1]{\F{\Sg^{\bullet}_{\bullet}}{\bar{#1}}}
\newcommand{\Dtcon}[1]{\F{\Sg^{\circ}_{\circ}}{\bar{#1}}}

\newcommand{\ACon}[1]{\FC{\R}{#1}}
\newcommand{\BCon}[1]{\FC{\R_{\bullet}}{\bar{#1}}}
\newcommand{\DCon}[1]{\FC{\R_{\circ}}{\bar{#1}}}
\newcommand{\AtCon}[1]{\FC{\Sg^\prime}{#1}}
\newcommand{\BtCon}[1]{\FC{\Sg^{\bullet}_{\circ}}{\bar{#1}}}
\newcommand{\CtCon}[1]{\FC{\Sg^{\bullet}_{\bullet}}{\bar{#1}}}
\newcommand{\DtCon}[1]{\FC{\Sg^{\circ}_{\circ}}{\bar{#1}}}

\newcommand{\Dconfat}[2]{\F{\R_{\circ}}{\bar{#1},\bar{#2}}}
\newcommand{\Btconfat}[2]{\F{\Sg_{\bullet}^{\circ}}{\bar{#1},\bar{#2}}}
\newcommand{\Dtconfat}[2]{\F{\Sg_{\circ}^{\circ}}{\bar{#1},\bar{#2}}}

\newcommand{\DConfat}[2]{\FC{\R_{\circ}}{\bar{#1},\bar{#2}}}
\newcommand{\BtConfat}[2]{\FC{\Sg_{\bullet}^{\circ}}{\bar{#1},\bar{#2}}}
\newcommand{\DtConfat}[2]{\FC{\Sg_{\circ}^{\circ}}{\bar{#1},\bar{#2}}}

%%%%%%%%%%%%%%%%%%%%%%%%%%%%%%%%%%%%%%%%%%%%%%%%%%%%%%%%%%%%%%%%%%%%%%%%%%%%%%%%%%%%%%%%%%%
%
%  paper formatting
%
%%%%%%%%%%%%%%%%%%%%%%%%%%%%%%%%%%%%%%%%%%%%%%%%%%%%%%%%%%%%%%%%%%%%%%%%%%%%%%%%%%%%%%%%%%%

\theoremstyle{plain}
\newtheorem{thm}{Theorem}[section]
\newtheorem{prop}[thm]{Proposition}
\newtheorem{cor}[thm]{Corollary}
\newtheorem{lem}[thm]{Lemma}
\newtheorem{conj}[thm]{Conjecture}

\theoremstyle{definition}
\newtheorem{defn}[thm]{Definition}
\newtheorem{exmp}[thm]{Example}

\theoremstyle{remark}
\newtheorem*{rem}{Remark}
\newtheorem*{hnote}{Historical Note}
\newtheorem*{nota}{Notation}
\newtheorem*{ack}{Acknowledgments}
\numberwithin{equation}{section}

\title {Particle configurations and Coxeter operads}

\subjclass[2000]{Primary 18D50, Secondary 14H10, 52B11}

\thanks{All authors were partially supported by NSF grant DMS-9820570.  Carr and Devadoss were also supported by NSF CARGO grant DMS-0310354.}

\author[Armstrong]{Suzanne M.\ Armstrong}
\address{S.\ Armstrong: Williams College, Williamstown, MA 01267}
\email{suzanne.m.armstrong@williams.edu}

\author[Carr]{Michael Carr}
\address{M.\ Carr: University of Michigan, Ann Arbor, MI 48109}
\email{mpcarr@umich.edu}

\author[Devadoss]{Satyan L.\ Devadoss}
\address{S.\ Devadoss: Williams College, Williamstown, MA 01267}
\email{satyan.devadoss@williams.edu}

\author[Engler]{Eric Engler}
\address{E.\ Engler: Williams College, Williamstown, MA 01267}
\email{eric.h.engler@williams.edu}

\author[Leininger]{Ananda Leininger}
\address{A.\ Leininger: MIT, Cambridge, MA 02139}
\email{anandal@mit.edu}

\author[Manapat]{Michael Manapat}
\address{M.\ Manapat: University of California, Berkeley, CA 94720}
\email{manapat@ocf.berkeley.edu}

\begin{abstract}
There exist natural generalizations of  the real moduli space of Riemann spheres based on manipulations of Coxeter complexes.  These novel spaces inherit a tiling by the graph-associahedra convex polytopes.  We obtain explicit configuration space models for the classical infinite families of finite and affine Weyl groups using particles on lines and circles.  A Fulton-MacPherson compactification of these spaces is described and this is used to define the Coxeter operad.  A complete classification of the building sets of these complexes is also given, along with a computation of their Euler characteristics.
\end{abstract}

\maketitle

%%%%%%%%%%%%%%%%%%%%%%%%%%%%%%%%%%%%%%%%%%%%%%%%%%%%%%%%%%%%%%%%%%%%%%%%%%%%%%%%%%%%%
%NOTE: The \baselineskip=15pt below is the minimum needed to make the paper look decent.
%      There are many superscripts {\tilde, \overline} which cause erratic
%      spaces between sentences. \baselineskip corrects these problems.  Also see the
%      \baselineskip=12pt by References.
%%%%%%%%%%%%%%%%%%%%%%%%%%%%%%%%%%%%%%%%%%%%%%%%%%%%%%%%%%%%%%%%%%%%%%%%%%%%%%%%%%%%%

\baselineskip=15pt

%%%%%%%%%%%%%%%%%%%%%%%%%%%%%%%%%%%%%%%%%%%%%%%%%%%%%%%%%%%%%%%%%%%%%%%%%%%%%%%%%%%%%
%
%                Motivation from Physics
%
%%%%%%%%%%%%%%%%%%%%%%%%%%%%%%%%%%%%%%%%%%%%%%%%%%%%%%%%%%%%%%%%%%%%%%%%%%%%%%%%%%%%%
\section{Motivation from Physics}
\subsection{}
A configuration space of $n$ ordered, distinct particles on a variety $V$ is
$$\F{V}{n} = V^n - \Delta, \ \ \ {\rm where} \ \Delta = \{(x_1, \ldots , x_n) \in V^n \suchthat \exists \ i,j, \ x_i = x_j \}.$$
Over the past decade, there has been an increased interest in the configuration space of $n$ labeled particles on the projective line. The focus is on a quotient of this space by $\PGLC$, the affine automorphisms on $\CP$.  The resulting variety is the moduli space of Riemann spheres with $n$ punctures ${\mathcal M}_{0, n} = \F{\Cx \Pj^1}{n}/\Pj \Gl_2(\Cx).$  There is a compactification \CM{n} of this space, a smooth variety of complex dimension $n-3$, coming from Geometric Invariant Theory \cite{git}.  The space \CM{n} plays a crucial role as a fundamental building block in the theory of Gromov-Witten invariants, also appearing in symplectic geometry and quantum cohomology \cite{km}.

Our work is motivated by the \emph{real} points \M{n} of this space, the  set of points fixed under complex conjugation.  These real moduli spaces have importance in their own right, appearing in areas such as $\zeta$-motives of Goncharov and Manin \cite{gm} and Lagrangian Floer theory of Fukaya \cite{foo}.  Indeed, \M{n} has even emerged in phylogenetic trees \cite{bhv} and networks \cite{lp}.  It was Kapranov \cite{kap1} who first noticed a relationship between \M{n} and the braid arrangement of hyperplanes, associated to the Coxeter group of type $A$:  Blow-ups of certain cells of the $A_n$ Coxeter complex yield a space homeomorphic to a double cover of \M{n+2}. This creates a natural tiling of \M{n} by associahedra, the combinatorics of which is discussed in \cite{dev3}.  Davis et.\ al have generalized this construction to all Coxeter groups, along with studying the fundamental groups of these blown-up spaces \cite{djs}.  Carr and Devadoss \cite{cd} looked at the inherent tiling of these spaces by the convex polytopes \emph{graph-associahedra}.

%%%%%%%%%%%%%%%%%%%%%%%%%%%%%%%%%%%%%%%%%%%%%%%%%%%%%%%%%%%%%%%%%%%%%%%%%%%%%%%%%%%%%
\subsection{}

We begin with elementary results and notation:  Section~\ref{s:complexes} provides the background of Coxeter groups and their associated Coxeter complexes and Section~\ref{s:config} constructs the appropriate configuration spaces of particles. Section~\ref{s:bracket} introduces the bracketing notation in order to visualize collisions in the configuration spaces, leading to viewing the hyperplanes of the Coxeter complexes in this new language.  It is this notation that provides a transparent understanding of several results in this paper.  In particular, it allows for a complete classification of the minimal building sets for the Coxeter complexes, along with their enumeration, as given in Tables~\ref{t:sphere} and~\ref{t:euclid}. It is interesting to note that some configuration structures behave quite classically, whereas others (based on \emph{thick} particles) are atypical.  

The heart of the paper begins in Section~\ref{s:compact} where the Fulton-MacPherson compactification \cite{fm} of these spaces is discussed and used to define the notion of a \emph{Coxeter operad} in Definition~\ref{d:coxoperad}, extending the mosaic operad of \M{n} \cite{dev1}.   Here, \emph{nested} bracketings are used to describe the structure of the compactified spaces, enabling us to describe how the chambers of these spaces glue together.  Section~\ref{s:tiling} ends with combinatorial results using the theory developed in \cite{cd}.  For instance, the Euler characteristics of these Coxeter moduli spaces are given, exploiting the tilings by graph-associahedra.

As of writing this paper, the importance of the operadic structure of \M{n} has been brought further to light.  For instance, Etingof et al.\ \cite{ehkr} have used the mosaic operad to compute the cohomology ring of \M{n}.  Recently, the Coxeter operad defined below appears in the work of E.\ Rains in computing the integral homology of these generalized Coxeter moduli spaces \cite{rai}.

\begin{ack}
We express our gratitude to Jim Stasheff and Vic Reiner for detailed comments on this work and to Alex Postnikov for working on the $f$-vectors of graph-associahedra.  We also thank the referee for wonderful insight and careful analysis, especially with regards to the operad module structures.
\end{ack}

%%%%%%%%%%%%%%%%%%%%%%%%%%%%%%%%%%%%%%%%%%%%%%%%%%%%%%%%%%%%%%%%%%%%%%%%%%%%%%%%%%%%%
%
%                Spherical and Euclidean Complexes
%
%%%%%%%%%%%%%%%%%%%%%%%%%%%%%%%%%%%%%%%%%%%%%%%%%%%%%%%%%%%%%%%%%%%%%%%%%%%%%%%%%%%%%
\section{Spherical and Euclidean Complexes}
\label{s:complexes}
\subsection{}

In order to provide the construction of Coxeter moduli space generalizations of \M{n}, as given in Section~\ref{s:compact}, we begin with standard facts and definitions about Coxeter systems.  Most of the material here can be found in Bourbaki \cite{bou}.

\begin{defn}
Given a finite set $S$, a \emph{Coxeter group} \ $W$ is given by the presentation
$$W \ = \ \langle \ s_i \in S \ \ | \ \ s_i^2 = 1, \ (s_i s_j)^{m_{ij}} = 1 \ \rangle \, ,$$
where $m_{ij} = m_{ji}$ and $2 \leq m_{ij} \leq \infty$. The pair $(W,S)$ is called a \emph{Coxeter system}.
\end{defn}

\noindent Associated to any Coxeter system $(W, S)$ is its \emph{Coxeter graph} $\Cox_W$:  Each node represents an element of $S$, where two nodes $s_i, s_j$ determine an edge if and only if $m_{ij} \geq 3$.
A Coxeter group is \emph{irreducible} if its Coxeter graph is connected and it is \emph{locally finite} if either $W$ is finite or each proper subset of $S$ generates a finite group.  A Coxeter group is \emph{simplicial} if it is irreducible and locally finite \cite{djs}. The classification of simplicial Coxeter groups and their graphs is well-known \cite[Chapter 6]{bou}.

\begin{figure}[h]
\includegraphics {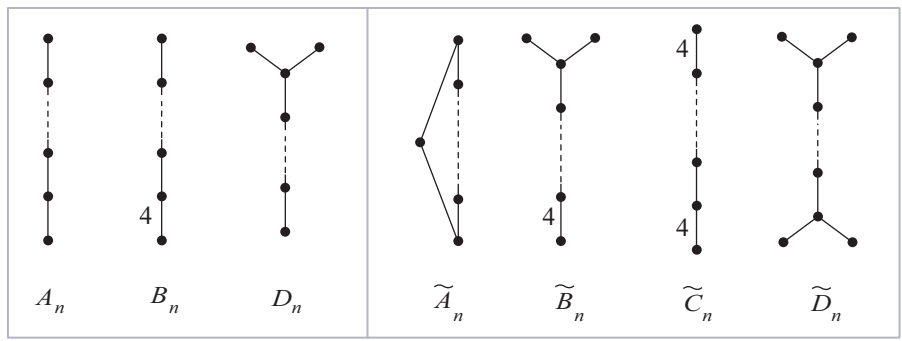}
\caption{Coxeter graphs of spherical and Euclidean groups.}
\label{f:coxdiags}
\end{figure}

We restrict our attention to infinite families of simplicial Coxeter groups which generalize to arbitrary number of generators.  This will mimic configuration spaces of an arbitrary number of particles since our motivation comes from \M{n}; it will also allow a well-defined construction of an operad.   There are only seven such types of Coxeter groups: three spherical ones and four Euclidean ones. Figure~\ref{f:coxdiags} shows the Coxeter graphs associated to the Coxeter groups of interest; we label the edge with its order for $m_{ij} > 3$.  The number of nodes of a graph is given by the subscript $n$ for the spherical groups, whereas the number of nodes is $n+1$ for the Euclidean case.

%%%%%%%%%%%%%%%%%%%%%%%%%%%%%%%%%%%%%%%%%%%%%%%%%%%%%%%%%%%%%%%%%%%%%%%%%%%%%%%%%%%%%
\subsection{}

Every spherical Coxeter group has an associated finite reflection group  realized by reflections across linear hyperplanes on a sphere.  Every conjugate of a generator $s_i$ acts on the sphere as a reflection in some hyperplane, dividing the sphere into simplicial chambers.  The sphere, along with its cellulation is the \emph{Coxeter complex} corresponding to $W$, denoted $\C{W}$.  The hyperplanes associated to each group given in Table~\ref{t:s} lie on the $(n-1)$ sphere.   The $W$-action on the chambers of $\C{W}$ is simply transitive, and thus we may associate an element of $W$ to each chamber. The number of chambers of $\C W$ comes from the order of the group.

\renewcommand{\arraystretch}{1.5}{
\begin{table}[h]
\begin{tabular}{crr}
$W$ & \hspace{1in} Hyperplanes &  \hspace{.8in} $\#$ Chambers \\ \hline
$A_n$ & $x_i = x_j$ &  \ $(n+1)!$ \\
$B_n$ &  $x_i = 0$, \ $x_i = \pm x_j$ & $2^{n} \ n!$\\
$D_n$ &  $x_i = \pm x_j$ & $2^{n-1} \ n!$
\end{tabular}
\bigskip
\caption{The spherical arrangements.}
\label{t:s}
\end{table}}
\vspace{-.2in}

We move from spherical geometry coming from linear hyperplanes to Euclidean geometry arising from affine hyperplanes.  Just as with the spherical case, each Euclidean Coxeter group has an associated Euclidean reflection group realized as reflections across affine hyperplanes in Euclidean space.  Again, we focus on the infinite families of such Euclidean Coxeter groups which are $\widetilde A_n, \widetilde B_n, \widetilde C_n$, and $\widetilde D_n$. The hyperplanes associated to each group, given in Table~\ref{t:t}, lie in $\R^n$.

We look at the quotient of the Euclidean space $\R^n$ by a group of translations, resulting in the $n$-torus $\T^n$.  This is done for three reasons:  First, the configuration space model is a more natural object after the quotient, resulting in particles on circles.  Second, it is the correct generalization of the affine type $A$ complex, as discussed in \cite{dev2}.  Third, and most importantly, it presents us with valid operad module structures as given in Section~\ref{s:compact}.

\renewcommand{\arraystretch}{1.5}{
\begin{table}[h]
\begin{tabular}{crr}
$W$ & \hspace{1.4 in} Hyperplanes ($k \in 2\Z$) & \hspace{.7in} $\#$ Chambers \\ \hline
$\widetilde A_n$ &  $x_i = x_j + k$ & $n!$\\
$\widetilde B_n$ &  $x_i = \pm x_j + k, \ x_i = 1 + k$ & $2^{n-1} \ n!$\\
$\widetilde C_n$ &  $x_i = \pm x_j + k, \ x_i = 1 + k, \ x_i = 0 + k$ & $2^{n} \ n!$\\
$\widetilde D_n$ &  $x_i = \pm x_j + k$ & $2^{n-2} \ n!$
\end{tabular}
\bigskip
\caption{The toroidal arrangements.}
\label{t:t}
\end{table}}
%\vspace{-.2in}

The translations for $\widetilde A_n$ are covered in \cite[Section 2.3]{dev2}.  For the remaining cases, we choose a slightly non-standard collection of hyperplanes in order for the associated configuration spaces to be more canonical.  This has the benefit of identifying the most ubiquitous set of hyperplanes $\{x_i = \pm x_j + k\}$,
producing an arrangement that is familiar from the spherical cases. We refer to the quotient of the complex $\C W$ of Euclidean type as the \emph{toroidal Coxeter complex}, denoted $\T\C{W}$. 

\begin{exmp}
Figure~\ref{f:2cpx}(a) is $\cas{3}$, the $2$-sphere with hyperplane markings, and part (b) is the 2-sphere $\cbs{3}$. Figure~\ref{f:2cpx}(c) shows the hyperplanes of $\widetilde C_2$ in $\R^2$, whereas (d) shows the hyperplanes for $\widetilde A_2$. Part (e) is the cellulation of the toroidal complex $\cat{2}$.

\begin{figure}[h]
\includegraphics[width=\textwidth]{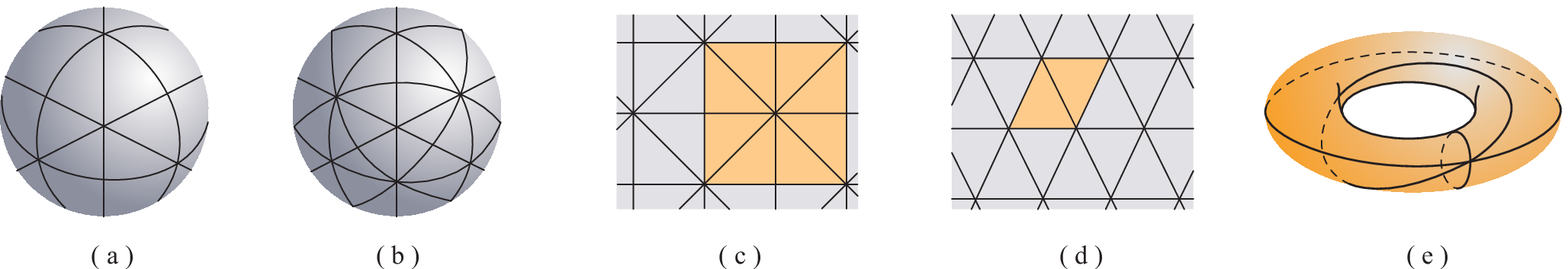}
\caption{Coxeter complexes (a) $\cas{3}$, (b) $\cbs{3}$, (c) $\C {\widetilde C_2}$, (d) $\C {\widetilde A_2}$ and (e) $\cat{2}$.}
\label{f:2cpx}
\end{figure}
\end{exmp}

%%%%%%%%%%%%%%%%%%%%%%%%%%%%%%%%%%%%%%%%%%%%%%%%%%%%%%%%%%%%%%%%%%%%%%%%%%%%%%%%%%%%%
%
%                Configuration Spaces
%
%%%%%%%%%%%%%%%%%%%%%%%%%%%%%%%%%%%%%%%%%%%%%%%%%%%%%%%%%%%%%%%%%%%%%%%%%%%%%%%%%%%%%
\section{Configuration Spaces}
\label{s:config}
\subsection{}

We now give an explicit configuration space analog to each Coxeter complex above. These appear as (quotients of) configuration spaces of particles on the line $\R$ and the circle $\Sg$.  The arguments used for the constructions below are elementary, immediately following from the hyperplane arrangements of the reflection groups.  However, as shown in Section~\ref{s:compact}, the configuration space model we provide will enable us to elegantly capture the blow-ups of these Coxeter complexes.

\begin{defn}
Let $\Acon{n} = \R^n - \{ (x_{1}, x_{2}, ..., x_{n}) \in \mathbb{R}^{n} \suchthat \exists \; i,j, \ x_{i}= x_{j} \}$ be the configuration space of $n$ labeled particles on the real line $\R$.  A generic point in $\F{\R}{5}$ is $ \ x_ 1 \ < \ x_2 \ < \ x_3 \ < \ x_4 \ < \ x_5 \ ,$
which we notate (without labels) as \ \includegraphics{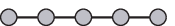}\ .
\end{defn}

\begin{defn}
Let $\Bcon{n} = \R^n - \{ (x_{1}, x_{2}, ..., x_{n}) \in \mathbb{R}^{n} \suchthat \exists \; i,j, \ x_{i}= \pm x_{j} \mbox{ or } x_{i}=0 \}$ be the space of $n$ pairs of \emph{symmetric} labeled particles (denoted $\bar{n}$) across the origin.  A point in $\Bcon{3}$ is $ \ -x_3 \ < \ -x_2 \ < \ -x_1 \ < \ 0 \ < \ x_1 \ < \ x_2 \ < \ x_3 \ ,$ which is depicted without labels as \includegraphics{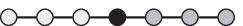}\ , where the black particle is fixed at the origin.
\end{defn}

\begin{defn}
Let $\Dcon{n} = \R^n - \{(x_{1}, x_{2}, ..., x_{n}) \in \mathbb{R}^{n} \suchthat  \exists \; i,j, \ x_{i}= \pm x_{j}\}$ be the space of $\bar n$ pairs of symmetric labeled particles across the origin, where the particle $x_i$ and its symmetric partner $-x_i$ are both allowed to occupy the origin. A point in $\Dcon{3}$ is
\begin{equation}
-x_3 \ < \ -x_2 \ < \ -x_1, x_1 \ < \ x_2 \ < \ x_3 \,
\label{e:points}
\end{equation}
drawn \ \includegraphics{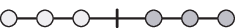} \ without labels.  Notice the mark at the origin where there is no fixed particle:  The point $ \ -x_3 \ < \ -x_2 \ < \ x_1, -x_1 \ < \ x_2 \ < \ x_3 \ $ drawn as \includegraphics{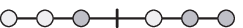} lies in the same chamber of $\Dcon{3}$ as Eq.\eqref{e:points}.
\end{defn}

Let $\Aff(\R)$ be the group of affine transformations of $\R$ generated by translating and positive scaling. The action of $\Aff(\R)$ on $\Acon{n}$ translates the leftmost of the $n$ particles in $\R$ to $-1$ and the rightmost is scaled to $1$.  If we allow the particles in $\Acon{n} / \Aff(\R)$ to \emph{collide} (coincide with each other), the resulting space is denoted $C_n \langle \R \rangle$.  In a sense, this includes the hyperplanes which were removed back into $\R^n$.  The space $C_n \langle \R \rangle$ is sometimes referred to as the \emph{naive compactification} of $\Acon{n} / \Aff(\R)$.\footnote{The space $C_n \langle \R \rangle$  can also be thought of as the \emph{closure} of $\Acon{n} / \Aff(\R)$, though the space in which this closure is taken is non-trivial.  An excellent treatment of these ideas in a general context is given by Sinha \cite[Section 3]{sin}.}

\begin{prop}
$C_n \langle \R \rangle$ has the same cellulation as $\cas{n-1}$.
\end{prop}

\noindent A detailed proof of this is given in \cite[Section 4]{dev3}.  Roughly, quotienting by translations of $\Aff(\R)$ removes the inessential component of the arrangement and scaling results in restricting to the sphere $\cas{n-1}$.  This proposition can be extended to the other spherical Coxeter complexes.  Let $\Aff(\bar \R)$ be the transformations of $\R$ generated simply by positive scalings:  The action of $\Aff(\bar \R)$ scales the (symmetric) particles farthest from the origin to unit distance.  Let $C_{\bar n} \langle \R_{\bullet} \rangle$ and $C_{\bar n} \langle \R_{\circ} \rangle$ denote spaces where particles of $\Bcon{n} / \Aff(\bar \R)$ and $\Dcon{n} / \Aff(\bar \R)$ have collided, respectively.

\begin{prop}
$C_{\bar n} \langle \R_{\bullet} \rangle$ and $C_{\bar n} \langle \R_{\circ} \rangle$ have the same cellulation as $\cbs{n}$ and $\cds{n}$ respectively.
\end{prop}

\begin{proof}
From the above definition above, it is clear $\Bcon{n}$ and $\Dcon{n}$ are complements of the hyperplanes given in Table~\ref{t:s}. Thus any collision of particles in the configuration space maps to a point on the hyperplanes defined by the associated finite reflection group.  Quotienting by $\Aff(\bar \R)$  allows choosing a particular representative for each fiber.  Specifically, for the fiber containing $(x_1, x_2, ..., x_n)$, choose
$$(x_1, x_2, ..., x_n)/\sqrt{x_1^2+x_2^2+...+x_n^2},$$
giving a map onto the unit sphere in $\R^n$.  The  cellulation of the sphere by these hyperplanes yields the desired Coxeter complex.
\end{proof}

%%%%%%%%%%%%%%%%%%%%%%%%%%%%%%%%%%%%%%%%%%%%%%%%%%%%%%%%%%%%%%%%%%%%%%%%%%%%%%%%%%%%%
\subsection{}

\begin{figure}[h]
\includegraphics{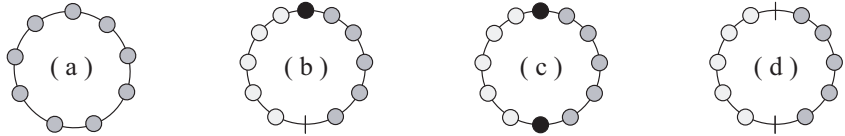}
\caption{Configurations of particles (a) without and (b) - (d) with symmetry.}
\label{f:ca-d}
\end{figure}

We move from the spherical to the affine (toroidal) complexes.  However, the interest now is on configurations of particles on the circle $\Sg$.  The group of rotations acts freely on $\F{\Sg}{n}$, and its quotient is denoted by $\Atcon{n}$; Figure~\ref{f:ca-d}(a) shows a point in $\Atcon{9}$ drawn without labels.

\begin{prop}
$\AtCon{n}$ has the same cellulation as $\cat{n-1}$.
\end{prop}

\noindent A proof of this is given in \cite[Section 3]{dev2}.  A similar construction is produced below for the other three toroidal Coxeter complexes.  Our focus now is on the circle $\Sg$ with the vertical line through its center as its axis of symmetry, where the two diametrically opposite points on the axis are labeled $0$ and $1$.  The space of interest is the configuration space of pairs of \emph{symmetric} labeled particles (again denoted $\bar{n}$) across this symmetric axis of the circle.

\begin{defn}
Let $\Btcon{n} = \T^n - \{ (x_{1}, x_{2}, ..., x_{n}) \in \T^{n} \suchthat \exists \: i,j, \ x_{i}= \pm x_{j}  \mbox{ or } x_{i}=1\}$ be the space of $n$ pairs of symmetric labeled particles on $\Sg$ with a fixed particle at $1$.  Figure~\ref{f:ca-d}(b) shows a point in $\Btcon{5}$.
\end{defn}

\begin{defn}
Let $\Ctcon{n} = \T^n -  \{ (x_{1}, x_{2}, ..., x_{n}) \in \T^{n} \suchthat \exists \: i,j, \ x_{i}= \pm x_{j} \mbox{ or } x_{i}=1 \mbox{ or } x_{i}=0 \}$ be the space of $n$ pairs of symmetric labeled particles on $\Sg$ with a fixed particle at $0$ and $1$. Figure~\ref{f:ca-d}(c) shows a point in $\Ctcon{5}$.
\end{defn}

\begin{defn}
Let $\Dtcon{n} = \T^n -  \{ (x_{1}, x_{2}, ..., x_{n}) \in \T^{n} \suchthat \exists \: i,j, \ x_{i}= \pm x_{j}\}$ be the space of $n$ pairs of symmetric labeled particles on $\Sg$ with no fixed particles. Figure~\ref{f:ca-d}(d) shows a point in $\Dtcon{5}$.
\end{defn}

\begin{prop}
$\BtCon{n}$, $\CtCon{n}$ and $\DtCon{n}$ have the same cellulation as $\cbt{n}$, $\cct{n}$ and $\cdt{n}$ respectively.
\end{prop}

\begin{proof}
This is a direct consequence of the definitions of the configuration spaces, of the toroidal complexes, and their corresponding hyperplane arrangements given in Table~\ref{t:t}.
\end{proof}

%%%%%%%%%%%%%%%%%%%%%%%%%%%%%%%%%%%%%%%%%%%%%%%%%%%%%%%%%%%%%%%%%%%%%%%%%%%%%%%%%%%%%
\subsection{}
\label{ss:strans}

The group of reflections $W$ across the respective hyperplanes acts on the configuration space by permuting particles.  The Coxeter group $A_n$ (the symmetric group)  is generated by transpositions $s_{ij}$ which interchange the $i$-th and $j$-th particle.  The Coxeter group $B_n$ of \emph{signed permutations} is generated by $s_{ij}$ along with reflections $r_1, \ldots, r_n$, where $r_i$ changes the sign of the $i$-th particle.  Note that $B_n$ is isomorphic to $\Z_2^n \,\rtimes\, \Sg_n$.
The Coxeter group $D_n$ is classically represented as the group of \emph{even} signed permutations.  Alternatively, $D_n$ is isomorphic to the group $\Z_2^{n-1} \,\rtimes\, \Sg_n$, generated by transpositions $s_{ij}$ along with reflections $r_2, \ldots, r_n$.  The element $r_1$, which is present in $B_n$ but not in $D_n$, corresponds to the reflection of the particle and its inverse that is closest to the origin.

Let $\sg(W)$ denote the group acting \emph{simply transitively} on the configuration spaces above.  As mentioned above, for the spherical Coxeter groups, $\sg(W)$ is isomorphic to $W$.  However, the action of the affine groups is only transitive on the toroidal complexes.  The simplest way to compute $\sg(W)$ for the toroidal cases is from observing the diagrams given in Figure~\ref{f:ca-d}.  Cutting the circle along a fixed point and ``laying it flat'' gives us the appropriate groups.  Since a particle in $\Atcon{n}$ is fixed by the group of rotations, then $\sg(\widetilde A_n)$ is isomorphic to $A_{n-1}$.  Similarly, $\sg(\widetilde B_n)$ is isomorphic to $D_n$ and $\sg(\widetilde C_n)$ is isomorphic to $B_n$.  The group $\sg(\widetilde D_n)$ is isomorphic to $\Z_2^{n-2} \,\rtimes\, \Sg_n$, generated by transpositions $s_{ij}$ along with reflections $r_2, \ldots, r_{n-1}$.  The elements $r_1$ and $r_n$, which are not present in $\sg(\widetilde D_n)$, correspond to the reflections of the particles and their inverses that are closest to the centrally symmetric axis.

%%%%%%%%%%%%%%%%%%%%%%%%%%%%%%%%%%%%%%%%%%%%%%%%%%%%%%%%%%%%%%%%%%%%%%%%%%%%%%%%%%%%%
%
%                Bracketings and Hyperplanes
%
%%%%%%%%%%%%%%%%%%%%%%%%%%%%%%%%%%%%%%%%%%%%%%%%%%%%%%%%%%%%%%%%%%%%%%%%%%%%%%%%%%%%%
\section{Bracketings and Hyperplanes}
\label{s:bracket}
\subsection{}

We introduce the bracket notation in order to visualize collisions in the configuration spaces, leading to a transparent understanding of our results below.  In particular, Proposition~\ref{p:bracket} produces a complete classification of the minimal building sets for the Coxeter complexes, along with their enumeration, as given in Tables~\ref{t:sphere} and~\ref{t:euclid}.

A \emph{bracket} is drawn around adjacent particles on a configuration space diagram representing the collision of the included particles.   A $k$-\emph{bracketing} of a diagram is a set of $k$ brackets representing multiple independent particle collisions.  For example, the configuration
\begin{equation}
-x_4 \ = \ -x_3 \ < \ -x_2 \ < \ -x_1 \ = \ 0 \ = \ x_1 \ < \ x_2 \ <  \ x_3 \ = \ x_4
\label{e:br}
\end{equation}
in $\BCon{5}$ corresponds to the bracketing \includegraphics{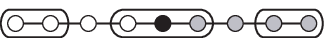}.  Each bracket on a configuration space diagram with symmetric particles will actually consist of two symmetric brackets, one on each side of the origin, with this symmetric pair counting as only one bracket. If this set includes the origin, we draw one symmetric bracket around the origin, which again counts as one bracket.  Thus Eq.\eqref{e:br} is a $2$-bracketing of its diagram.

Let $\al$ be an intersection of hyperplanes. We say that hyperplanes $h_i$ \emph{cellulate} $\al$ to mean the intersections $h_i \cap \al$ decompose $\al$ into cells. Denote by $\Hs{\al}$ the set of all hyperplanes that contain $\al$, called the stabilizing hyperplanes of $\al$.  If reflections in these hyperplanes generate a finite reflection group, it is called the \emph{stabilizer} of $\al$. Note that in a simplicial Coxeter complex, the stabilizer exists for all intersections of hyperplanes.\footnote{By abuse of terminology, we also refer to the set $\Hs{\al}$ as the stabilizer of $\al$.}

We define the \emph{support} of a bracketing to be the configuration space associated to the bracketing diagram.  That is, it is the subspace (of the configuration space) in which particles that share a bracket have collided.  However, a set of collisions in a configuration space defines an intersection of hyperplanes.  So, alternatively, the support of a bracketing is the smallest intersection of hyperplanes associated to the bracketing.  The following is immediate:

\begin{lem}
If $\al$ is the support of a bracketing $G$, then for every pair of particles $x_i$ and $x_j$ that share a bracket in $G$, the hyperplane defined by $x_i=x_j$ is in $\Hs \al$.
\label{l:stab}
\end{lem}

\begin{figure}[h]
\includegraphics[width=\textwidth] {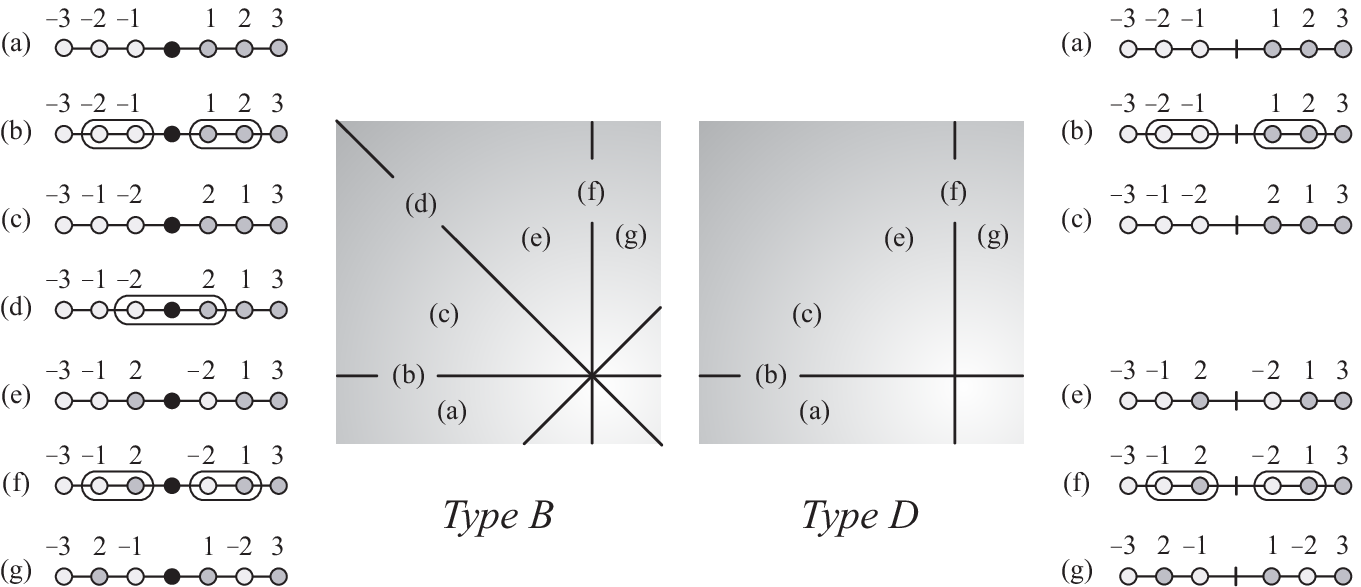}
\caption{Local regions of $\cbs{3}$ and $\cds{3}$.}
\label{f:typebd}
\end{figure}

\begin{exmp}
Figure~\ref{f:typebd} shows part of the two-dimensional complexes $\cbs{3}$ and $\cds{3}$, one with and one without a fixed particle at the axis of symmetry. As we move through the chambers, going from (a) through (g), a representative of each configuration is shown.  Notice that since there is no fixed particle at the axis of symmetry for type $D$, there is no meaningful bracketing of the symmetric particles closest to the axis; they may pass each other freely \emph{without} collision.
\end{exmp}

%%%%%%%%%%%%%%%%%%%%%%%%%%%%%%%%%%%%%%%%%%%%%%%%%%%%%%%%%%%%%%%%%%%%%%%%%%%%%%%%%%%%%
\subsection{}

There are natural composition maps on bracketed diagrams.  These form the basic operations of our operads defined in Section~\ref{s:compact}.

\begin{defn}
There are three types of compositions:
\begin{enumerate}
\item Let $H$ be a diagram of $m$ particles of $\Acon{k}$ or $\Atcon{k}$, with one particle labeled $i$.  Let $G$ be a diagram of $\Acon{k}$. The \emph{composition} $H \circ_i G$ is the diagram of $m+k-1$ particles where the particle $i$ is replaced by a bracket containing $G$.

\item Let $H$ be a diagram of $m$ paired particles of a configuration space, with one particle labeled $i$ and its mirror image labeled $-i$.  Let $G$ be a diagram of $\Acon{k}$. The \emph{composition} $H \circ_i G$ is the diagram of $m+k-1$ paired particles, where the particle $i$ is replaced by a bracket containing $G$ and its pair $-i$ is replaced by a bracket containing the mirror image of $G$ (left side of Figure~\ref{f:bracketcomp}).

\begin{figure}[h]
\includegraphics {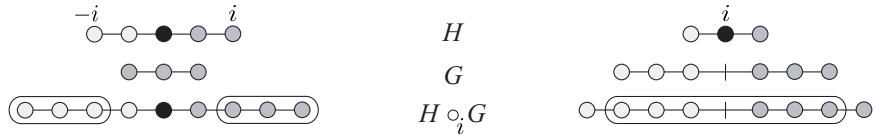}
\caption{Composition operations on bracketings.}
\label{f:bracketcomp}
\end{figure}

\item Let $H$ be a diagram of $m$ paired particles of $\Bcon{m}$, $\Btcon{m}$ or $\Ctcon{m}$, with a fixed particle labeled $i$.  Let $G$ be a diagram of either $\Bcon{k}$ or $\Dcon{k}$. The \emph{composition} $H \circ_i G$ is the diagram of $m+k$ paired particles where the fixed particle $i$ is replaced by a bracket containing $G$ (right side of Figure~\ref{f:bracketcomp}).
\end{enumerate}
\label{d:comp}
\end{defn}

\noindent Indeed any $k$-bracketing $G$ can be represented as
\begin{equation}
G \ = \ H \circ_{i_1} G_1 \circ_{i_2}  \cdots \circ_{i_k}G_k,
\label{e:bracket}
\end{equation}
where the base $H$ and the $G_i$'s are diagrams without brackets and each $i_j$ is a particle in $H$.

\begin{prop} \label{p:bracket}
Let $G$ be a $k$-bracketing as defined in Eq.\eqref{e:bracket}.  Moreover, let $\al$ be the support of $G$ and let $\al_i$ be the support of $G_i$.  Then the stabilizer of $\al$ is the product of the stabilizers of $\al_i$ and the cellulation of $\al$ is determined by $H$.
\end{prop}

\begin{proof}
The product structure of the stabilizer of $\al$ is a consequence of Lemmas~\ref{l:stab}.  The cellulation of $\al$ is determined by the hyperplanes of $\al$, which are simply the intersections of hyperplanes of $\C W$ with $\al$.  If the particles $x_i$ and $x_j$ share a bracket in $\al$, then $x_k$ cannot collide with $x_j$ in $\al$ without also colliding with $x_i$. Geometrically, this property implies that the two hyperplanes $x_i=x_k$ and $x_j=x_k$ have the same intersection with $\al$. Similarly, since $x_i$ and $x_j$ collide in all of $\al$, the hyperplane $x_i=x_j$ plays no role in the cellulation of $\al$.  These two facts allow us to treat $x_i$ and $x_j$ as a single particle in $G$ without changing the hyperplane arrangement. Repeating this process for all particles that share a bracket gives the desired result.
\end{proof}

%%%%%%%%%%%%%%%%%%%%%%%%%%%%%%%%%%%%%%%%%%%%%%%%%%%%%%%%%%%%%%%%%%%%%%%%%%%%%%%%%%%%%
\subsection{}

It is easy to check that in most cases the cellulations of subspaces (intersections of hyperplanes) in Coxeter complexes are indeed other (smaller dimensional) Coxeter complexes.  There are, however, three instances where this is not so, appearing as subspaces of the Coxeter complexes  $\cds{n}$, $\cbt{n}$ and $\cdt{n}$.  In particular, they have cellulations combinatorially equivalent to Coxeter complexes with \emph{additional} hyperplanes.  We define these three atypical complexes below:

\begin{defn}  The complexes of interest are:

\noindent 1. \ \ Let $\csds{n}{m}$ be $\cds{n}$ with $m$ additional hyperplanes $\{x_i = 0 \suchthat 1 \le i \le m \}$.

\noindent 2. \ \ Let $\csbt{n}{m}$ be $\cbt{n}$ with $m$ additional hyperplanes $\{x_i = 0 \suchthat 1 \le i \le m \}$.

\noindent 3. \ \ Let $\csdt{n}{m}$ be $\cdt{n}$ with $2m$ additional hyperplanes $\{x_i = 0, 1 \suchthat 1 \le i \le m \}$.

\label{d:atypical}
\end{defn}

\begin{figure}[h]
\includegraphics {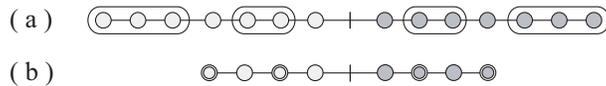}
\caption{Bracketing and thick particles.}
\label{f:thick}
\end{figure}

The configuration space model provides intuition into how these cases arise naturally.  Note how these are all complexes with associated configuration spaces on $\R_{\circ}$, $\Sg^{\bullet}_{\circ}$ and $\Sg^{\circ}_{\circ}$, where not all points along the axis of symmetry have fixed particles.  The subspaces of these configuration spaces are those where some particles have collided. In these subspaces, sets of collided particles may be considered in aggregate as a new type of particle, called a \emph{thick} particle. Figure~\ref{f:thick}(a) shows a bracketing and (b) its representation with thick particles.  In general, thick particles allow us to represent any number of coincident particles by a single particle.

Recall that particles were defined such that they could occupy the same point as their inverse; that is, they do not form a collision with their inverse.  Unlike (standard) particles, a thick particle and its inverse may not occupy the same point without collision.   The reason comes from the hyperplane equations:  In the subspace where $x_i$ and $x_j$ have collided, the hyperplane $x_i = -x_j$ represents the same configurations as $x_i = -x_i \ ( = 0)$.  Thus, the $m$ additional hyperplanes added to the complex correspond to the $m$ thick particles in their configuration spaces.  Then the  diagram of Figure~\ref{f:thick}(b) is an element of $\csds{4}{2}$, sitting as a subspace of $\cds{7}$ in Figure~\ref{f:thick}(a).

\begin{rem}
In the case of non-paired particles, the distinction between standard and thick particles is irrelevant, since no particle has an inverse to collide with.  They are also inconsequential in configuration spaces that include a fixed particle wherever particles may meet their inverses.
\end{rem}

%%%%%%%%%%%%%%%%%%%%%%%%%%%%%%%%%%%%%%%%%%%%%%%%%%%%%%%%%%%%%%%%%%%%%%%%%%%%%%%%%%%%%
%
%                Compactifications and Operads
%
%%%%%%%%%%%%%%%%%%%%%%%%%%%%%%%%%%%%%%%%%%%%%%%%%%%%%%%%%%%%%%%%%%%%%%%%%%%%%%%%%%%%%
\section{Compactifications and Operads}
\label{s:compact}
\subsection{}

Compactifying a configuration space $\F{V}{n}$ enables the points on $V$ to collide and a {\em system} is introduced to record the {\em directions} points arrive at the collision.  In the work of Fulton and MacPherson~\cite{fm}, this method is brought to rigor in the algebro-geometric context.\footnote{Axelrod and Singer \cite{as} look at this compactification from a perspective of spherical blow-ups on real manifolds.} In \cite[Section 4]{dp}, De Concini and Procesi show that the \emph{minimal blow-ups} of the Coxeter complexes $\C W$ are equivalent to the Fulton-MacPherson compactifications of their corresponding configuration spaces.

In order to describe these compactified Coxeter moduli spaces, we begin with definitions.  The collection of hyperplanes $\{x_i = 0 \ | \ i = 1, \ldots, n\}$ of $\R^n$ generates the \emph{coordinate} arrangement.  A crossing of hyperplanes is \emph{normal}\, if it is locally isomorphic to a coordinate arrangement.
A construction which transforms any crossing into a normal crossing involves the algebro-geometric concept of a blow-up; a standard reference is \cite{har}.

\begin{defn}
The \emph{blow-up} of a space $V$ along a codimension $k$ intersection $\al$ of hyperplanes is the closure of $\{(x,f(x))\suchthat x\in V\}$ in $V \times \Pj^{k-1}$. That is, we replace $\al$ with the projective sphere bundle associated to the normal bundle of $\al$.
\end{defn}

A general collection of blow-ups is usually noncommutative in nature; in other words, the order in which spaces are blown up is important.  For a given arrangement, De Concini and Procesi \cite[Section 3]{dp} establish the existence and uniqueness of a \emph{minimal building set}, a collection of subspaces for which blow-ups commute for a given dimension, and for which every crossing in the resulting space is normal. For a Coxeter complex $\C{W}$, we denote the minimal building set by $\Min{W}$.

\begin{defn}
The \emph{Coxeter moduli space} $\Cm{W}$ is the minimal blow-up of $\C{W}$, obtained by blowing up along elements of $\Min{W}$ in \emph{increasing} order of dimension.
\end{defn}

\begin{rem}
Kapranov showed the minimal blow-ups of $\C A_n$ yield a space homeomorphic to a double cover of the moduli space \M{n+2} \cite[Proposition 4.8]{kap1}.  Thus, the Fulton-MacPherson compactifications of our configuration spaces yield generalizations of this moduli space.
\end{rem}

\begin{figure}[h]
\includegraphics [width=\textwidth]{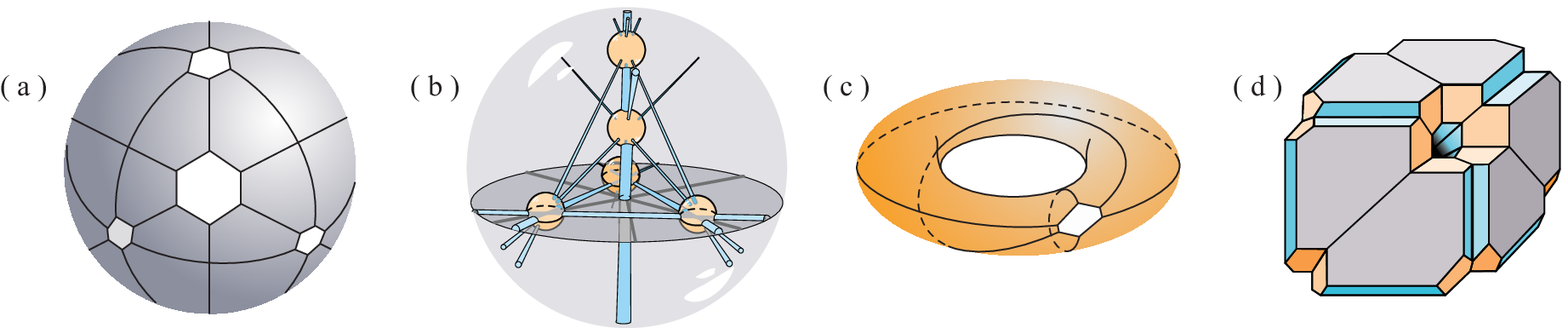}
\caption{Coxeter moduli spaces (a) $\Cm{A_3}$, (b) $\Cm{A_4}$, (c) $\Cm{\widetilde A_2}$ and (d) $\Cm{\widetilde A_3}$.}
\label{f:type2blow}
\end{figure}

\begin{exmp}
Figure \ref{f:type2blow}(a)  shows the blow-ups of the sphere $\C{A_3}$ of Figure~\ref{f:2cpx}(a) at nonnormal crossings.  Each blown up point has become a hexagon with antipodal identification and the resulting manifold is the Coxeter moduli space $\Cm{A_3}$, homeomorphic to the eight-fold connected sum of $\R\Pj^2$ with itself.\footnote{A projective version of this diagram is first found in a different context by Brahana and Coble in $1926$~\cite[Section 1]{bc} relating to possibilities of maps with twelve five-sided countries.}  Part (b) shows \M{6} coming from the iterated blow-ups of $\C{A_4}$.  Figure~\ref{f:type2blow}(c) shows the minimal blow-up of $\C{\widetilde A_2}$ of Figure~\ref{f:2cpx}(e).  Finally, part (d) is the cube with opposite facial identifications yielding the three-torus of $\Cm{\widetilde A_3}$.
\end{exmp}

The relationship between the set $\Min{W}$ and the group $W$ is given by the concept of reducibility.  For intersections of hyperplanes $\al,\be,$ and $\ga$, the collection of hyperplanes $\Hs \al$ is \emph{reducible} if $\Hs \al$ is a disjoint union $\Hs \be \sqcup \Hs \ga$, where $\al = \be \cap \ga$.

\begin{lem} \cite[Section 3]{djs}
Let $\al$ be an intersection of hyperplanes of $\C W$.   Then $\Hs{\al}$ is irreducible if and only if $\al$ is in $\Min{W}$.
\label{l:irred}
\end{lem}

\noindent This lemma can be rewritten in the language of bracketings:

\begin{lem}
Let $\al$ be an intersection of hyperplanes of $\C W$. Then $\al$ is in $\Min{W}$ if and only if $\al$ is the support of a $1$-bracketing.
\label{l:conneq}
\end{lem}

\begin{proof}
If $\al$ is the support of a $1$-bracketing $G$, then $\Hs{\al}$ is determined  by Lemma~\ref{l:stab}. Thus:
\begin{enumerate}
\item If $G$ contains a fixed particle, then $\Hs{\al} \cong \Hy B_k$.

\item If $G$ contains a particle and its inverse but no fixed particle, then $\Hs{\al} \cong \Hy D_k$.

\item If $G$ does not contain a particle and its inverse, then $\Hs{\al} \cong \Hy A_k$.
\end{enumerate}
\noindent All three of these hyperplane arrangements are irreducible, so $\al$ is in $\Min{W}$.

Conversely, let $\al$ be the support of a $k$-bracketing $G = H \circ_{i_1} G_1 \circ_{i_2}  \cdots \circ_{i_k}G_k$.  Let $\be$ be the support of $G_1$ and $\ga$ be the product of the configuration spaces diagramed by $G_2, \cdots, G_k$. By the definition of reducibility, $\Hs{\al}= \Hs{\be} \sqcup \Hs{\ga}$, and thus $\al$ is not in $\Min{W}$ by Lemma \ref{l:irred}.
\end{proof}

\begin{rem}
Table~\ref{t:sphere} itemizes the collection of elements in $\Min{W}$ for the spherical cases and Table~\ref{t:euclid} for the Euclidean ones.  In the tables, $m$ represents the total number of thick particles and $r$ the number of thick particles in the bracket (stabilizer) of the atypical complexes.
\end{rem}

%%%%%%%%%%%%%%%%%%%%%%%%%%%%%%%%%%%%%%%%%%%%%%%%%%%%%%%%%%%%%%%%%%%%%%%%%%%%%%%%%%%%%
\subsection{}

As bracketings encoded collisions of particles in configuration spaces, it is \emph{nested} bracketings which encode the Fulton-MacPherson compactification of the configuration spaces \cite[Section 2]{dp}. The FM compactification allows collisions of particles whose description comes from the repulsive potential observed by quantum physics:  Pushing particles together creates a spherical bubble onto which the particles escape \cite{pw}.  In other words, as particles try to collide, the result is a new bubble fused to the old at the point of collision, where the collided particles are now on the new bubble.  The phenomena is dubbed as {\em bubbling}, with the resulting structure as a bubble-tree.   Indeed, the nested bracketings are exactly the one-dimensional analogues of bubble-trees \cite[Section 1]{dev1}.  

\begin{figure}[h]
\includegraphics[width=.8\textwidth]{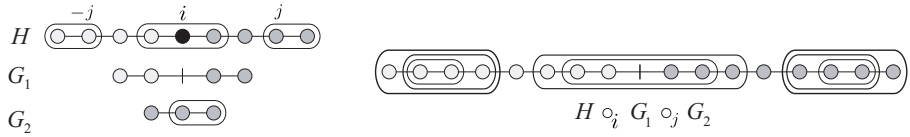}
\caption{Composition operations on nested bracketings.}
\label{f:nestbracket}
\end{figure}

Moreover, the codimension $k$ faces of a chamber of a compactified configuration space are the nested $k$-bracketings on the configuration space diagrams.  A diagram $G$ denoted
\begin{equation}
G \ = \ H \circ_{i_1} G_1 \circ_{i_2} \cdots \circ_{i_k} G_k,
\label{e:nestbrac}
\end{equation}
is a nested  $(k + m_0 + m_1 + \cdots + m_k)$-bracketing, with each $i_j$ a particle of $H$, where $H$ is an $m_0$-bracketing and where $G_i$ is a nested $m_i$-bracketing; see Figure~\ref{f:nestbracket}.  The composition maps are those in Definition~\ref{d:comp}.

After compactification, different orderings of particles in a bracket do not necessarily represent the same cell. A different action is necessary to describe the identification of diagrams.

\begin{defn}
The \emph{flip} $\sghat(G)$ action on an unbracketed diagram $G$ consists of the identity and the reflection (that reverses the order of the particles of $G$). On a nested bracketing diagram, the action of $\sghat(G)$ acts independently on each bracketed component.
\end{defn}

\begin{thm}
\label{t:flip}
Let $G$ be a nested bracketing where $G \ = \ H \circ_{i_1} G_1 \circ_{i_2} \cdots \circ_{i_k} G_k$. \
All bracketings in the image of \ $G$ under $\sghat(G_1)  \times \cdots \times \sghat(G_k)$ represent the same cell.
\end{thm}

\begin{proof}
By definition, blow-ups introduce a \emph{projective} bundle around each subspace in the minimal building set.  The analog in configuration spaces is an identification across each bracket:  Flipping the positions of the particles in the bracket is defined to represent the same configuration.  Thus the permutations that represent the same set of configurations as $G$ are exactly the images of $G$ under $\sghat(G_1)  \times \cdots \times \sghat(G_k)$.
\end{proof}

\begin{rem}
This theorem gives us a gluing rule between the faces of two chambers.  In other words,
two nested $k$-bracketings $G_1$ and $G_2$ of a diagram (representing codimension $k$-faces) are identified if $G_2$ can be obtained from $G_1$ by flipping some of the brackets of $G_2$.
\end{rem}

\begin{figure}[h]
\includegraphics[width=\textwidth] {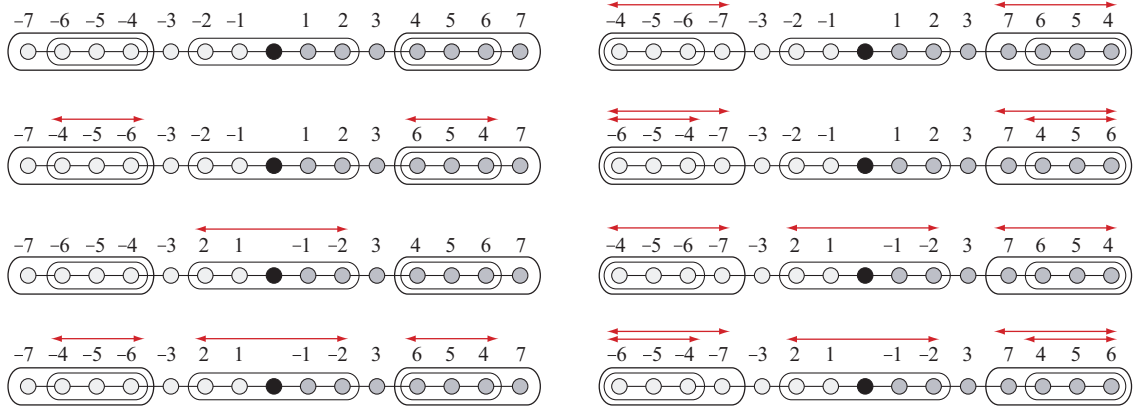}
\caption{Flips on nested bracketings.}
\label{f:flip}
\end{figure}

\noindent Figure \ref{f:flip} shows the permutations that preserve the cell represented by a particular configuration space diagram.  Note that reflections commute with each other, and thus they generate a group which is isomorphic to $(\Z/2\Z)^3$.

%%%%%%%%%%%%%%%%%%%%%%%%%%%%%%%%%%%%%%%%%%%%%%%%%%%%%%%%%%%%%%%%%%%%%%%%%%%%%%%%%%%%%
\subsection{}

Classically, the notion of an operad was created for the study of iterated loop spaces \cite{sta}.  Since then, operads have been used as universal objects representing a wide range of algebraic concepts.  An \emph{operad} $\Op$ consists of a collection of objects $\{\Op(n) \; | \; n \in \mathbb {N} \}$ in a monoidal category endowed with certain extra structures.  Notably, $\Op(n)$ carries an action by the symmetric group of $n$ letters, and there are composition maps
$$\Op(n) \otimes \Op(k_1) \otimes \cdots \otimes \Op(k_n) \rightarrow \Op(k_1 + \cdots + k_n)$$
which must be associative, unital, and equivariant; see \cite[Chapter 1]{mss} for details.

\begin{figure}[h]
\includegraphics[width=.9\textwidth]{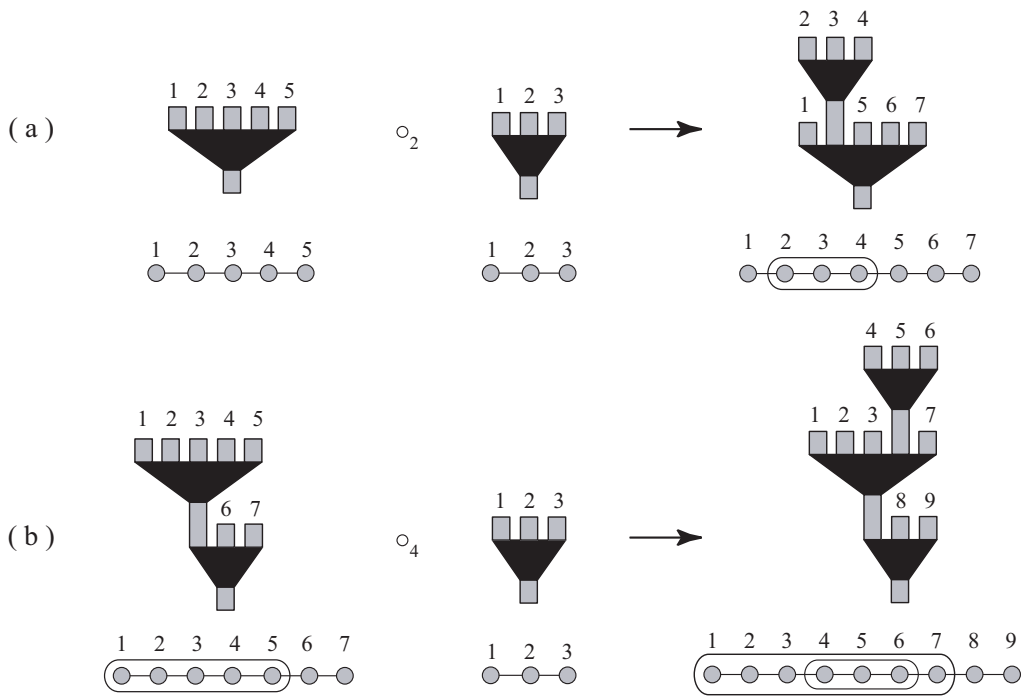}
\caption{Examples of composition maps of an operad, along with dual figures related to (a) bracketing and (b) nested bracketing.}
\label{f:operad}
\end{figure}

One can view $\Op(n)$ as objects consisting of $n$-ary operations, which yield an output given $n$ inputs.  We will be concerned mostly with operads in the context of topological spaces, where the objects $\Op(n)$ will be equivalence classes of geometric objects.  Classically, these objects can be pictured as in Figure~\ref{f:operad}.  The composition $\Op(i) \circ_k \Op(j)$ is obtained by grafting the output of $\Op(j)$ to the $k$-th input of $\Op(i)$.  The symmetric group acts by permuting the labeling of the inputs.

There are several variants and extensions of the operad definition above.  
A classic example is the \emph{non-symmetric} version, which removes all references to the symmetric group in the definition above.  
A \emph{right module} $\Op_R$ over an operad $\Op$ is a collection of objects $\{\Op_R(n) \; | \; n \in \mathbb {N} \}$ with a set of composition maps
$$\Op_R(n) \otimes \Op(k_1) \otimes \cdots \otimes \Op(k_n) \rightarrow \Op_R(k_1 + \cdots + k_n).$$
A \emph{bi-colored} operad is a multicategory with two objects \cite[Section 2]{mar}:  Intuitively, each of the inputs and output is given one of two colors.  An element $\Op(j)$ can be grafted into an input of $\Op(i)$ if and only if the colors of the corresponding input and output match.   
Our version (the naming of which is credited to J.\ Stasheff) replaces the symmetric group in the classical definition with appropriate Coxeter groups instead.  Let $\tot$ be the collection $\{A, B, D, \widetilde A, \widetilde B, \widetilde C, \widetilde D\}$ of Coxeter groups from Figure~\ref{f:coxdiags}. 

\begin{defn}
For $W \in \tot$, let $\Op_W(n,k)$  be the collection of configuration spaces of nested $k$-bracketings with $n$ particles associated to $\Cm{W}$.  If $\Opu$ is a subset of $\tot$, then let $\Op_\Opu = \bigcup_{W \in \, \Opu} \Op_{W}.$  Then, the \emph{Coxeter operad} is defined for each pair $\Opu, \Opv$ of  (possibly empty) subsets of $\tot$ for which the collection
$$\Op_\Opu(n_H, k_H) \otimes \Op_\Opv(n_1, k_1) \otimes \cdots \otimes \Op_\Opv(n_m, k_m) \rightarrow \Op_\Opu(n_*, k_*)$$
of composition maps exist, where $n_* = n_H -m + \sum n_i$ and $k_* = k_H + m + \sum k_i$.
\label{d:coxoperad}
\end{defn}

\noindent We note several structures that appear based on this definition, giving a partial list: 

\bigskip
\paragraph{\emph{Classic Operads:}\ }
When $\Opu = \Opv = \{A\}$, the Coxeter operad becomes the $A_\infty$ operad structure \cite{sta} of the associahedron.  The classic symmetric group acting on $\Op_A$ is exactly the $A_n$ Coxeter group.  Examples of this are seen in Figure~\ref{f:operad}.

\bigskip
\paragraph{\emph{Right Modules:}\ }
When $\Opv = \{A\}$ and $\Opu = \{W\}$, for any $W \in \tot \setminus A$, we have a right module $\Op_W$ over the operad  $\Op_A$.  The composition map (based on non-nested bracketings) is described in  Definition~\ref{d:comp} (2) for centrally symmetric spaces.  Moreover, Definition~\ref{d:comp} (1) provides the composition for $W=\widetilde A$, resulting in the cyclohedral structure given in \cite[Section 4.4]{mss}.

\begin{figure}[h]
\includegraphics[width=.9\textwidth]{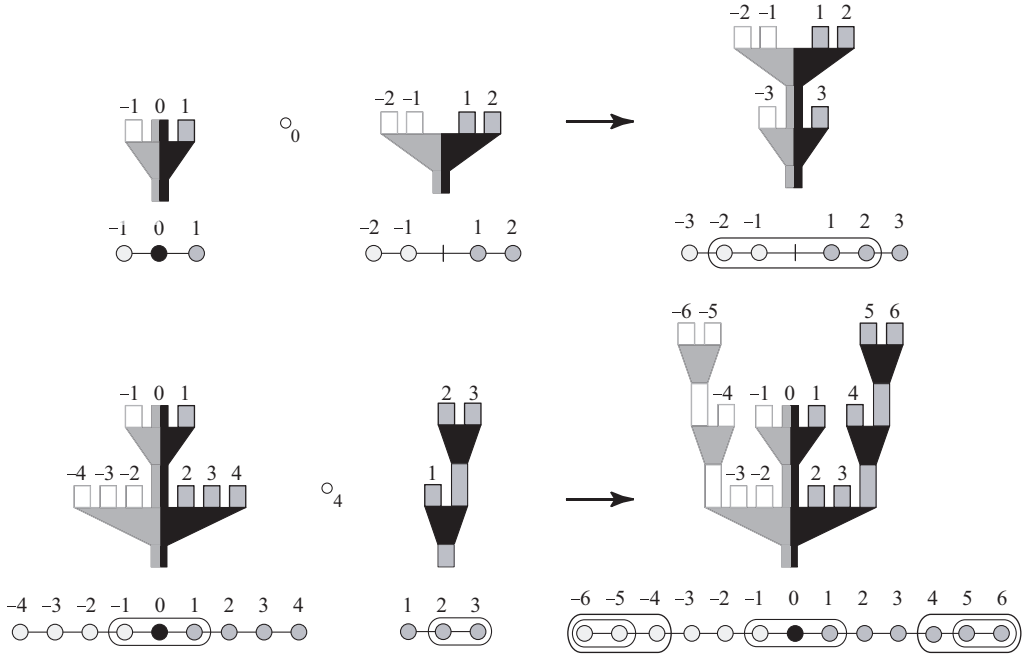}
\caption{Examples of composition maps of the Coxeter operad, along with dual \emph{tree} figures.}
\label{f:trees}
\end{figure}

\bigskip
\paragraph{\emph{Bi-Colored Operads:}\ } 
When $\Opu = \{B\}$ and $\Opv = \{A, B, D\}$, we obtain a bi-colored operad.  The two colors come from the centrally symmetric black particle and the ordinary free particles.  The map for $D \in \Opv$ is diagrammed in the top part of Figure~\ref{f:trees}, whereas the map for $A \in \Opv$ is given in the bottom of the figure.   This composition map is described in Definition~\ref{d:comp} (3).

\bigskip
\paragraph{\emph{Bi-Colored Right Modules:}\ }
The examples where $\Opu = \{\widetilde B, \widetilde C, \widetilde D\}$ and $\Opv = \{A, B, D\}$ result in the composition of bracketings on lines being glued onto bracketings on circles.   This results in a bi-colored right module over a bi-colored operad.  There are several options here which work for the different subsets of $\Opu$ and $\Opv$ given.  For example, $\Op_{\widetilde D}$ admits an action of $\Op_{A}$ and $\Op_{\{\widetilde B, \widetilde C\}}$ admits an action of $\Op_{\{B, D\}}$.  From an operadic viewpoint, the affine Coxeter groups are analogous to the spherical ones shown, but an \emph{unrooted tree}, rather than a rooted one, is used.

\begin{rem}
There exists a generalization of the classical operad to the \emph{braid} operad, defined by Fiedorowicz, with the braid group playing the role of the symmetric group \cite[Section 3]{zig}.  Since the braid group is the Artin group of type $A$, it seems plausible that the Coxeter operads above can be extended to their corresponding Artin groups, yielding analogs to the braid operads.
\end{rem}

%%%%%%%%%%%%%%%%%%%%%%%%%%%%%%%%%%%%%%%%%%%%%%%%%%%%%%%%%%%%%%%%%%%%%%%%%%%%%%%%%%%%%
%
%                Tiling by Graph-Associahedra
%
%%%%%%%%%%%%%%%%%%%%%%%%%%%%%%%%%%%%%%%%%%%%%%%%%%%%%%%%%%%%%%%%%%%%%%%%%%%%%%%%%%%%%
\section{Tiling by graph-associahedra}
\label{s:tiling}
\subsection{}

This section uses the theory of graph-associahedra developed in \cite{cd}.   As associahedra are related to the $A_\infty$ operad, we show that graph-associahedra capture the structure of the Coxeter operad.  Moreover, we provide combinatorial and enumerative results about the Coxeter moduli spaces.  Notably, the Euler characteristics of these spaces are given.

\begin{defn}
Let $\Cox$ be a graph.  A \emph{tube} is a proper nonempty set of nodes of $\Cox$ whose induced graph is a proper, connected subgraph of $\Cox$.  There are three ways that two tubes $t_1$ and $t_2$ may interact on the graph.
\begin{enumerate}
\item Tubes are \emph{nested} if  $t_1 \subset t_2$.
\item Tubes \emph{intersect} if $t_1 \cap t_2 \neq \emptyset$ and $t_1 \not\subset t_2$ and $t_2 \not\subset t_1$.
\item Tubes are \emph{adjacent} if $t_1 \cap t_2 = \emptyset$ and $t_1 \bigcup t_2$ is a tube in $\Cox$.
\end{enumerate}
Tubes are \emph{compatible} if they do not intersect and they are not adjacent.  A \emph{tubing} $T$ of $\Cox$ is a set of tubes of $\Cox$ such that every pair of tubes in $T$ is compatible.  A \emph{$k$-tubing} is a tubing with $k$ tubes.
\end{defn}

\begin{thm} \cite[Section 3]{cd}
For a graph $\Cox$ with $n$ nodes, the \emph{graph-associahedron} $\PG{}$ is the convex polytope of dimension $n-1$ whose face poset is isomorphic to set of valid tubings of \,$\Cox$, ordered such that $T \prec T'$ if $T$ is obtained from $T'$ by adding tubes.
\end{thm}

\begin{figure}[h]
\resizebox{.9\textwidth}{!}{\includegraphics {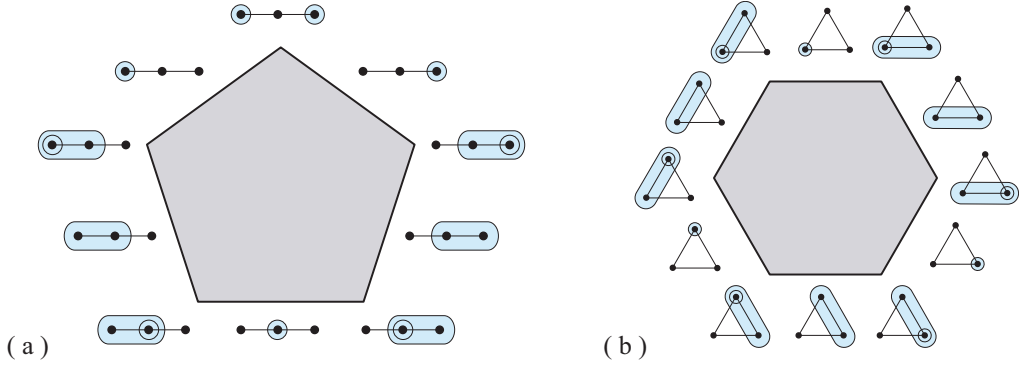}}
\caption{Graph-associahedra with a (a) path and (b) cycle as underlying graphs.}
\label{f:KWtubes}
\end{figure}

\noindent Figure~\ref{f:KWtubes} shows two examples of graph-associahedra, having underlying graphs as paths and cycles, respectively, with three nodes.  These turn out to be the two-dimensional associahedron \cite{sta} and cyclohedron \cite{bt} polytopes.

\begin{thm} \cite[Section 4]{cd} \label{t:tile}
Let $W$ be a simplicial Coxeter group and $\Cox_W$ be its associated Coxeter graph.  Then the $W$-action on $\C W_{\#}$ has a fundamental domain combinatorially isomorphic to $\PG{W}$.
\end{thm}

\begin{nota}
We write $\Pol W$ instead of $\PG{W}$ when context makes it clear.
\end{nota}

\begin{thm}
The tiling of the Coxeter moduli spaces are given as follows:
\begin{enumerate}
\item $\PA{n}$ (the associahedron) tiles $\Cas{n}$, $\Cbs{n}$ and $\Cct{n-1}$.
\item $\PAt{n}$ (the cyclohedron) tiles $\Cat{n}$.
\item $\PD{n}$ tiles $\Cds{n}$ and $\Cbt{n-1}$.
\item $\PDt{n}$ tiles $\Cdt{n}$.
\end{enumerate}
\label{t:tiling}
\end{thm}

\begin{proof}
For a given graph $\Cox$, the polytope $\PG{n}$  depends only on the adjacency of nodes, not the label on the edges.
\end{proof}

\begin{figure}[h]
\resizebox{.8\textwidth}{!}{\includegraphics {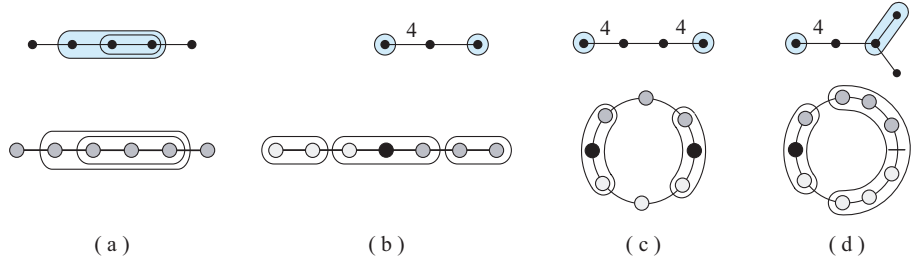}}
\caption{Examples of the bijection between tubings and nested bracketings.}
\label{f:tubebracket}
\end{figure}

\begin{rem}
There is a natural bijection from the set of all bracketings of a configuration space diagram to the set of all tubings of the associated Coxeter diagram.  The bijection is such that two brackets intersect if and only if their images intersect or are adjacent as tubes. Thus  the face poset of tubings is isomorphic to the face poset of bracketings, where $k$ brackets correspond to a codimension $k$ face.  Figure \ref{f:tubebracket} shows some  examples of this bijection.
\end{rem}

%%%%%%%%%%%%%%%%%%%%%%%%%%%%%%%%%%%%%%%%%%%%%%%%%%%%%%%%%%%%%%%%%%%%%%%%%%%%%%%%%%%%%

\subsection{}
We analyze the structure of these tiling polyhedra $\Pol {W}$.  For a given tube $t$ and a graph $\Cox$, let $\Cox_{t}$ denote the induced subgraph on the graph $\Cox$.  By abuse of notation, we sometimes refer to $\Cox_{t}$ as a tube.

\begin{defn}
Given a graph $\Cox$ and a tube $t$, construct a new graph $\Coxrec{t}$ called the \emph{reconnected complement}: If $V$ is the set of nodes of $\Cox$, then $V-t$ is the set of nodes of $\Coxrec{t}$.  There is an edge between nodes $a$ and $b$ in $\Coxrec{t}$ if either $\{a,b\}$ or $\{a,b\} \cup t$ is connected in $\Cox$.
\end{defn}

\begin{thm} \cite[Section 3]{cd}
\label{t:codim1}
The facets of \,$\PG{}$ correspond to the set of \,$1$-tubings on $\Cox$.  In particular, the facet associated to a $1$-tubing $\{t\}$ is equivalent to $\PG{t} \times \PGre{t}$.
\end{thm}

\noindent The facets of $\Pol W$ are of the form $\PG{} \times \PGre{}$, which can be found by simple inspection.  Using induction on each term of the product produces the following results:

\begin{cor} \cite{sta}
The faces of $\PA{}$ are of the form $\PA{} \times \cdots \times \PA{}$.
\end{cor}

\begin{cor} \cite{sta2}
The faces of $\PAt{}$ are of the form $\PAt{} \times \PA{} \times \cdots \times \PA{}$.
\end{cor}

Before moving on to the other tiling polytopes, we need to look at some special graphs which appear as reconnected complements. They are displayed in the Figure~\ref{f:specialdiags} below, the subscript $n$ denoting the number of vertices.  Note that the polytope $\PFp$ is the $3$-dimensional \emph{permutohedron}.

\begin{figure}[h]
\includegraphics {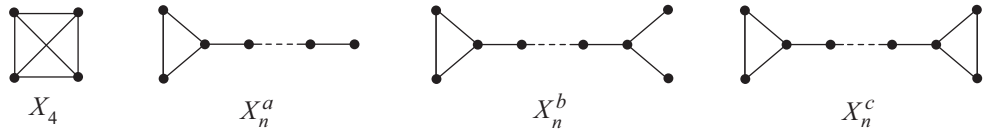}
\caption{Special graphs appearing as reconnected complements.}
\label{f:specialdiags}
\end{figure}

\begin{cor}
The faces of $\PD{}$ are of the form
\begin{enumerate}
\item $\PA{} \times \cdots \times \PA{}$ \item $\PD{} \times \PA{} \times \cdots \times \PA{}$ \item $\Pf{} \times \PA{} \times \cdots \times \PA{}.$
\end{enumerate}
\end{cor}

\begin{figure}[h]
\includegraphics[width=.9\textwidth] {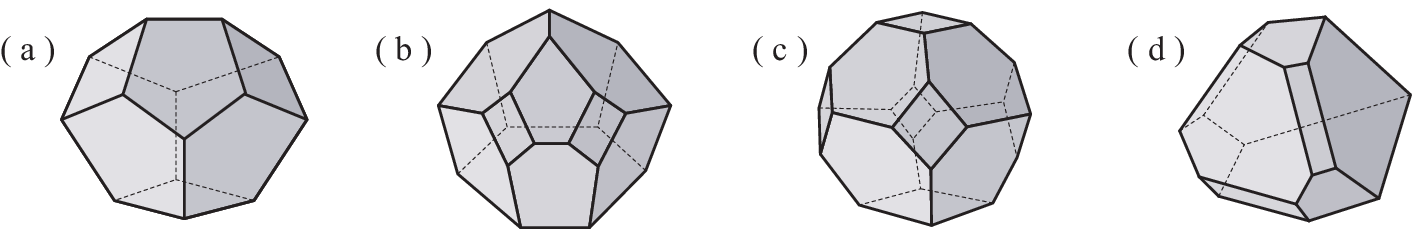}
\caption{The $3$-dimensional (a) associahedron $\PA{4}$, (b) cyclohedron $\PAt{3}$, (c) permutohedron $\PFp$ and (d) $\PD{4}$.}
\label{f:k5w4p4d4}
\end{figure}

\begin{exmp}
Figure~\ref{f:k5w4p4d4} illustrates four different polyhedra.  The first three are well-known objects:  (a) the associahedron $\PA{4}$, (b) cyclohedron $\PAt{4}$ and (c) permutohedron $\PFp$.  The last one (d) is $\PD{4}$ with six pentagons $\PA{3}$, three squares $\PD{2} \times \PA{2}$ and one hexagon $\Pf{3}$ for facets.
\end{exmp}

\begin{cor}
The faces of $\PDt{}$ are of the form
\renewcommand{\arraystretch}{1.33}{
$$\begin{array}{ll}
(1) \ \ \PA{} \times \cdots \times \PA{} &  (6) \ \ \PD{} \times \PD{} \times \PA{} \times \cdots \times \PA{} \\
(2) \ \ \PD{} \times \PA{} \times \cdots \times \PA{} & (7) \ \ \Pft{} \times \PA{} \times \cdots \times \PA{} \\
(3) \ \ \Pf{} \times \PA{} \times \cdots \times \PA{} & (8) \ \ \Pdf{} \times \PA{} \times \cdots \times \PA{} \\
(4) \ \ \Pf{} \times \Pf{} \times \PA{} \times \cdots \times \PA{} & (9) \ \ \PDt{} \times \PA{} \times \cdots \times \PA{} \\
(5) \ \ \PD{} \times \Pf{} \times \PA{} \times \cdots \times \PA{} & (10) \ \ \PFp \times \PA{} \times \cdots \times \PA{}.
\end{array}$$}
\end{cor}

%%%%%%%%%%%%%%%%%%%%%%%%%%%%%%%%%%%%%%%%%%%%%%%%%%%%%%%%%%%%%%%%%%%%%%%%%%%%%%%%%%%%%

\subsection{}
We compute the Euler characteristics of the Coxeter moduli spaces.  From Theorem~\ref{t:codim1}, we see that the number of codimension $k$ faces of the polytope $\Pol W$ cellulating $\C W_{\#}$ is precisely the number of $k$-tubings of the associated Coxeter diagram $\Cox_W$.

\begin{thm}
\label{t:euler}
Let $f_k(\Pol W)$ be the number of $k$-dimensional faces of $\Pol W$, and let $g$ be the number of chambers in the spherical or toroidal Coxeter complex $\C W$.  If\, $\dim(\Pol W) = n$, then
$$\chi(\,\Cm{W}\,) = \sum^{n}_{k = 0} \ (-1)^k \ \frac{g \cdot f_k}{2^{n - k}}.$$
\end{thm}

\begin{proof}
In order to count the number of $k$-dimensional faces in the space $\Cm{W}$, take the number of total chambers $g$ and multiply it by the number of $k$-dimensional faces $f_k(\Pol W)$ for each tile $\Pol W$.  The amount of overcounting is simply how different chambers of the complex meet at each $k$-dimensional face of a tile.  Since all the crossings in $\Cm{W}$ are normal, each $k$-dimensional face is identified with $2^{n-k}$ copies.
\end{proof}

\begin{rem}
The number of vertices in $\PA{}$ is the well-known Catalan number \cite[Section 6.5]{s2}.  The faces $f_k(W)$ of $\Pol W$ provide natural generalizations; see \cite{pos} for further exposition.
\end{rem}

Theorem~\ref{t:tiling} shows only four types of graph-associahedra tiling the Coxeter moduli spaces: $\PA{n}$ (the associahedron), $\PAt{n}$ (the cyclohedron), $\PD{n}$ and $\PDt{n}$.  The enumeration of the faces of the associahedra $\PA{n}$ is a classic result of A.\ Cayley \cite{cay}, obtained by just counting the number of $n$-gons with $k$ non-intersecting diagonals:
$$f_k(\PA{n}) = \frac{1}{n+1} \binom{n-1}{k} \binom{2n-k}{n}.$$
The enumeration of the face poset of the cylohedron $\PAt{n}$ comes from Simion \cite[Section 3]{sim}:
$$f_k(\PAt{n}) = \binom{n}{k} \binom{2n-k}{n}.$$
In a recent paper, Postnikov provides a recursive formula for the generating function of the numbers $f_k$ \cite[Theorem 7.11]{pos}.   Using this, a closed formulas for the graph-associahedra of types $D_n$ and $\widetilde D_n$ can be found; see \cite[Section 12]{pos2}.  We thank A.\ Postnikov for sharing the following result:

\begin{prop}
The face poset enumerations of types $D_n$ and $\widetilde D_n$ are
\begin{equation}
\label{e:postd}
f_k(\PD{n}) \ = \ 2 f_k(\PA{n}) - 2 f_k(\PA{n-1}) - f_k(\PA{n-2}) - f_{k-1}(\PA{n-1}) - f_{k-1}(\PA{n-2})
\end{equation}
\begin{equation}
\label{e:postdt}
\begin{split}
f_k(\PDt{n}) \ = \ & 4 f_k(\PA{n+1}) - 8 f_k(\PA{n}) -4 f_{k-1}(\PA{n}) + f_{k-2}(\PA{n-1}) \\
& \quad + 4 f_k(\PA{n-2}) + 6 f_{k-1}(\PA{n-2}) + 2 f_{k-2} (\PA{n-2}) \\
& \quad + f_k(\PA{n-3}) + 2 f_{k-1}(\PA{n-3}) + f_{k-2}(\PA{n-3}).
\end{split}
\end{equation}
\end{prop}

\begin{thm}
The Euler characteristics of the spherical blown-up Coxeter complexes are as follows: When $n$ is even, the values are zero; when $n=2m+1$ is odd,
\begin{eqnarray}
\label{e:as}
\chi(\,\Cas{n} \,) &=& (-1)^{m} \ 2 n \ ((n-2)!!)^2 \\
\label{e:bs}
\chi(\,\Cbs{n} \,) &=& 2^{n} \, \frac{1}{(n+1)} \ \chi(\,\Cas{n} \,) \\
\label{e:ds}
\chi(\,\Cds{n} \,) &=& 2^{n-3} \, \left[ \frac{8}{n+1} - \frac{1}{n-2} \right] \  \chi(\,\Cas{n} \,).
\end{eqnarray}
The Euler characteristics of the toroidal blown-up Coxeter complexes are as follows: When $n$ is odd, the values are zero; when $n=2m$ is even,
\begin{eqnarray}
\label{e:at}
\chi(\,\Cat{n} \,) &=& (-1)^{m} \ ((n-1)!!)^2 \\
\label{e:bt}
\chi(\,\Cbt{n} \,) &=& \frac{1}{2(n+1)} \ \chi(\,\Cds{n+1} \,) \\
\label{e:ct}
\chi(\,\Cct{n} \,) &=& 2^n \, \frac{1}{(n+2)(n+1)} \ \chi(\,\Cas{n+1} \,) \\
\label{e:dt}
\chi(\,\Cdt{n} \,) &=& 2^{n-6} \, \frac{1}{(n+1)} \, \left[ \frac{64}{n+2} - \frac{15}{n-1} \right] \  \chi(\,\Cas{n+1} \,).
\end{eqnarray}
\end{thm}

\begin{proof}
We use Theorem~\ref{t:euler} to obtain a summation, using the values $f_k$ given above along with the number of chambers provided by Tables~\ref{t:s} and \ref{t:t}.
%It is not hard to notice that alternating terms vanish.
The values for Eqs.\eqref{e:as} and \eqref{e:at} have been previously calculated in \cite[Section 3.2]{dev1} and \cite[Section 4.3]{sim} respectively.  Equations \eqref{e:bs}, \eqref{e:bt} and \eqref{e:ct} are consequences of Theorem~\ref{t:tiling}, where these spaces share the same tiling polytopes as previous calculations.
From Eq.~\eqref{e:postd} and Theorem~\ref{t:euler}, we obtain a linear combination
\begin{equation*}
\begin{split}
\chi(\,\Cds{n} \,) & = \frac{2^n}{(n+1)!} \ \chi(\,\Cas{n} \,) - \frac{2^{n-1}}{n!} \ \chi(\,\Cas{n-1} \,) - \frac{2^{n-3}}{(n-1)!} \ \chi(\,\Cas{n-2} \,) \\
& \quad + \frac{2^{n-1}}{n!} \ \chi(\,\Cas{n-1} \,) + \frac{2^{n-2}}{(n-1)!} \ \chi(\,\Cas{n-2} \,).
\end{split}
\end{equation*}
Algebraic manipulations result in Eq.\eqref{e:ds}.  Similar calculations using Eq.~\eqref{e:postdt} yield Eq.~\eqref{e:dt} after simplification.
\end{proof}

\begin{rem}
The reason there is a dimension shift between the spherical and toroidal cases is due to the convention of the affine case having $n+1$ nodes in its Coxeter graph, compared to $n$ nodes for the spherical.
\end{rem}

%%%%%%%%%%%%%%%%%%%%%%%%%%%%%%%%%%%%%%%%%%%%%%%%%%%%%%%%%%%%%%%%%%%%%%%%%%%%%%%%%%%%%
\subsection{}

The polytopes tiling the Coxeter moduli spaces are given by Theorem~\ref{t:tiling}.  We now discuss the tiling of the atypical complexes, given in Definition~\ref{d:atypical}, after minimal blow-ups.  As in other (compactified) configuration spaces, the chambers of these complexes correspond to orderings of the particles.  However, different orderings of particles may give different face posets to the chamber, since switching a standard and thick particle may change the valid bracketings of the diagram.  Specifically, near an axis of symmetry with no fixed particles, having thick particles allows more collisions and hence more brackets, than having standard particles.  It is here where the polytopes $\Pf{}, \Pft{},$ and $\Pdf{}$ based on Figure~\ref{f:specialdiags} appear.

Recall that the chambers of these complexes arise as subspaces of $\Cds{n}$, $\Cbt{n}$ or $\Cdt{n}$.  From Theorem~\ref{t:tiling}, these chambers must be \emph{faces} of either $\PD{}$ or $\PDt{}$.  Converting bracketings to tubings allows us to compute the face poset of the chamber using Theorem~\ref{t:codim1}.  The following is an example of this method.  Note how there is not simply one type of polytope tiling each blown-up atypical complex, as was the case with the Coxeter complexes.

\begin{figure}[h]
\includegraphics[width=.9\textwidth] {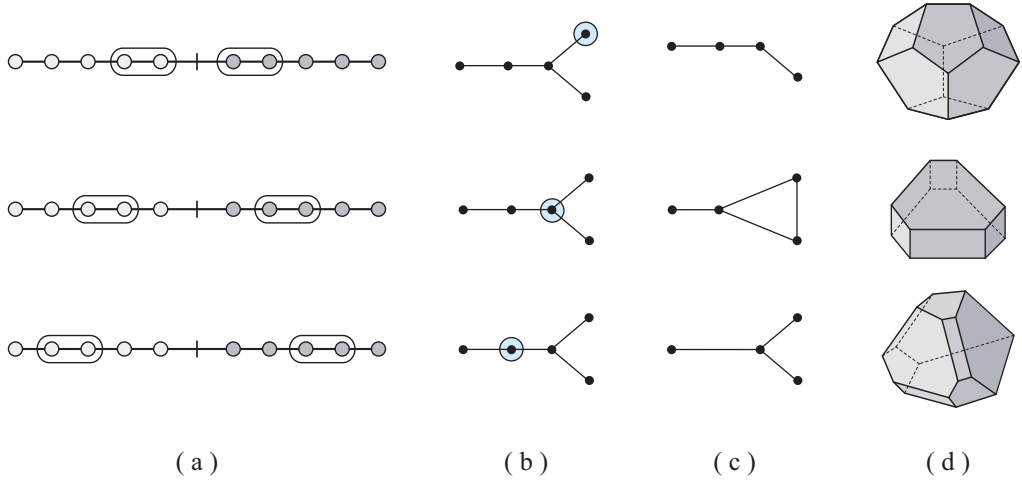}
\caption{(a) Configuration diagram brackets, (b) associated tubing, (c) reconnected complements and (d) fundamental chambers.}
\label{f:bracketrecon}
\end{figure}

\begin{exmp}
By taking different bracketings in configurations $\Dcon{5}$ associated to $\Cds{5}$, we can produce the configuration space diagrams of different chambers of $\Csds{4}{1}$, as in Figure~\ref{f:bracketrecon}.  By converting bracketings to tubings using the bijection, each facet of the chamber corresponds to the appropriate reconnected complement.
\end{exmp}

The gluing rules for the compactified configuration spaces applies to these atypical models as well.  The reflection action can change the ordering of standard and thick particles, and thus it encodes the manner in which polytopes of different types glue to tile the space.

\begin{figure}[h]
\includegraphics[width=\textwidth] {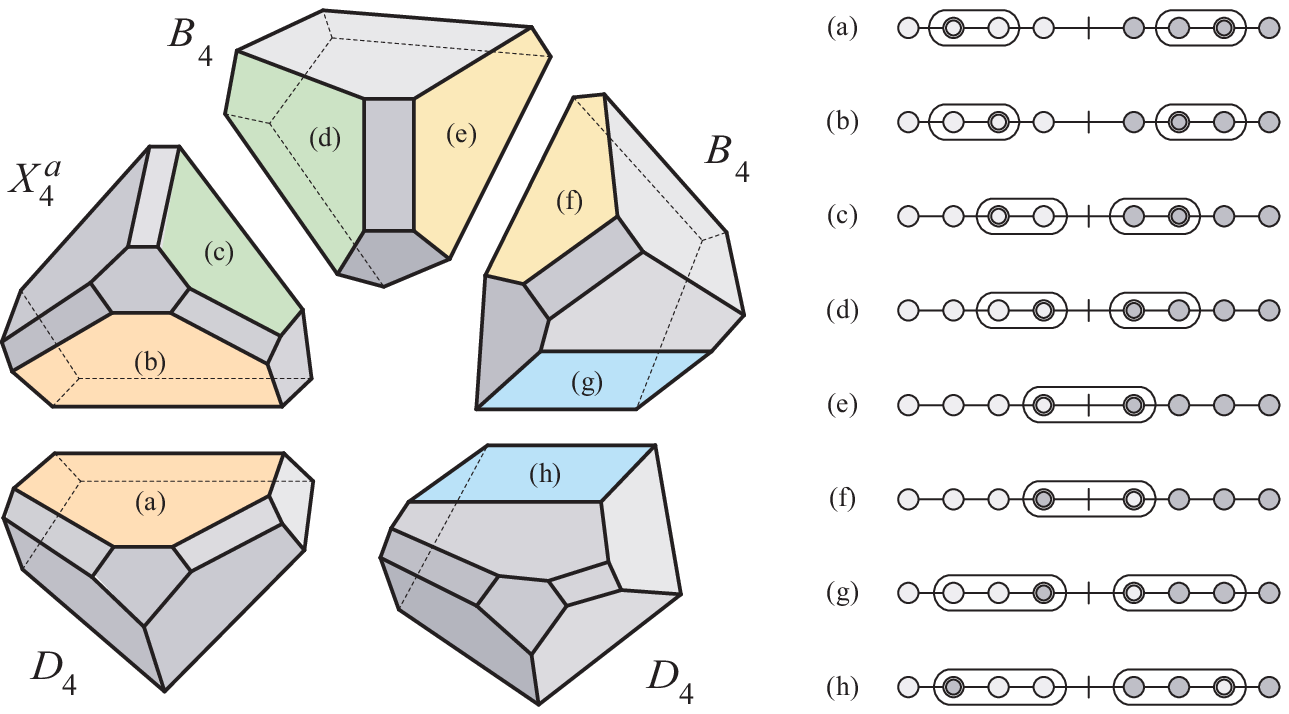}
\caption{Five adjacent chambers in $\Csds{4}{1}$.}
\label{f:dfbbd}
\end{figure}

\begin{exmp}
The illustration in Figure \ref{f:dfbbd} of five chambers of $\Csds{4}{1}$ shows how the configuration of three (standard) particles and 1 thick particle encodes gluing of faces among different types of chambers either across hyperplanes or antipodal maps. This is done by the flip action $\hat{\sigma}$ of Theorem~\ref{t:flip}.  We see the gluing of a face of $\PD{4}$ to a face of $\Pf{4}$ (a - b) and the corresponding labeled configuration spaces. Another face of this $\Pf{4}$ attaches to a face of $\PB{4}$ (c - d) which glues to the face of another $\PB{4}$ (e - f). This identification is across the hyperplane $x_{i}=0$, where $x_i$ is the label for the thick particle. Finally, this $\PB{4}$ glues to another chamber of type $\PD{4}$ (g - h) through an antipodal map.
\end{exmp}

%%%%%%%%%%%%%%%%%%%%%%%%%%%%%%%%%%%%%%%%%%%%%%%%%%%%%%%%%%%%%%%%%%%%%%%%%%%%%%%%%%%%%
%
%                  TABLES
%
%%%%%%%%%%%%%%%%%%%%%%%%%%%%%%%%%%%%%%%%%%%%%%%%%%%%%%%%%%%%%%%%%%%%%%%%%%%%%%%%%%%%%

\begin{table}[b]
\begin{tabular}{c | c c c c}
$\C W$ & Subspace & Stabilizer & Enumeration & Configuration \\ \hline \\

$\cas{n}$ &  $\cas{k+1}$ & $A_{n-k-1}$ &  $\binom{n+1}{n-k}$ & \resizebox{5cm}{!}{\includegraphics{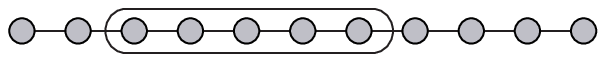}} \\ \\ \hline \\

        &  $\cbs{k+1}$ & $B_{n-k-1}$ &  $\binom{n}{n-k-1}$ & \includegraphics{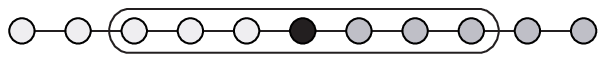}  \\
$\cbs{n}$ & & & &  \\
        &  $\cbs{k+1}$ & $A_{n-k-1}$ &  $2^{n-k-1} \binom{n}{n-k}$ & \includegraphics{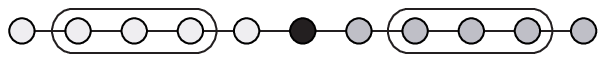}  \\ \\ \hline \\

        &  $\cbs{k+1}$ & $D_{n-k-1}$ &  $\binom{n}{n-k-1}$ & \includegraphics{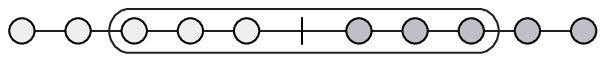}  \\
$\cds{n}$ & & & &  \\
        &  $\csds{k+1}{1}$ & $A_{n-k-1}$ &  $2^{n-k-1} \binom{n}{n-k}$ & \includegraphics{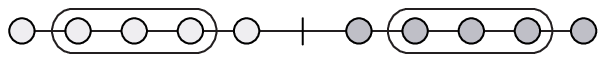}  \\ \\ \hline \\
        &  $\cbs{k+1}$ & $D_{{n-k-1},{r}}$ &  $\binom{m}{r} \binom{n-m}{n-k-r-1}$ & \includegraphics{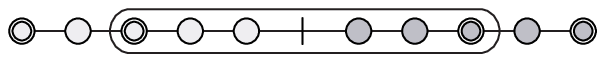} \\
$\csds{n}{m}$ & & & &  \\
        &  $\csds{k+1}{m-r+1}$ & $A_{n-k-1}$ &  $2^{n-k-1} \binom{n}{n-k}$ & \includegraphics{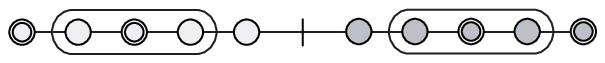}  \\ \\ \hline
\end{tabular}
\medskip
\caption{Minimal building sets of dimension $k$ for spherical complexes.}
\label{t:sphere}
\end{table}

%\begin{table}[h]
%\resizebox{\textwidth}{!}{\includegraphics{table1.eps}}
%\caption{Minimal building sets of dimension $k$ for spherical complexes.}
%\label{t:sphere}
%\end{table}

%%%%%%%%%%%%%%%%%%%%%%%%%%%%%%%%%%%%%%%%%%%%%%%%%%%%%%%%%%%%%%%%%%%%%%%%%%%%%%%%%%%%%

\begin{landscape}
\renewcommand{\arraystretch}{1}{
\begin{table}[h]
\begin{center}
\begin{tabular}{c | c | c c c | c c}
$\C W$ & $\cat{n}$ & &  $\cbt{n}$ & & \multicolumn{2}{|c}{$\cct{n}$} \\ \hline
& & & & & & \\
Subspace & $\cat{k+1}$ & $\csbt{k+1}{1}$ &  $\cbt{k+1}$ &  $\cct{k+1}$ & $\cct{k+1}$ & $\cct{k+1}$ \\
& & & & & & \\
Stabilizer & $A_{n-k}$ &  $A_{n-k-1}$  & $B_{n-k}$  & $D_{n-k}$ & $A_{n-k-1}$ & $B_{n-k}$ \\
& & & & & & \\
Enumeration  & $\binom{n+1}{n+1-k}$  &  $2^{n-k-1}\binom{n}{n-k}$ &  $\binom{n}{n-k}$ &  $\binom{n}{n-k}$  &
$2^{n-k-1}\binom{n}{n-k}$  &  $2 \binom{n}{n-k}$ \\
& & & & & & \\
\raisebox{.4in}{Configuration} & \includegraphics{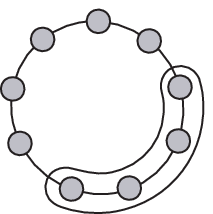} &  \includegraphics{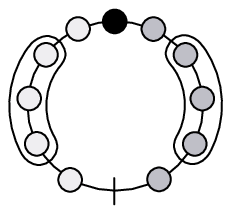} &
\includegraphics{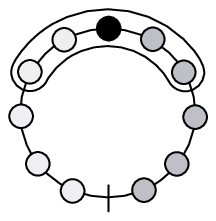} &  \includegraphics{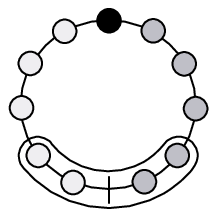} &  \includegraphics{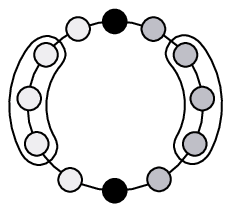} &
\includegraphics{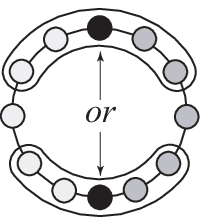} \\
\end{tabular}

\vspace{.5in}

\begin{tabular}{c | c c | c c c | c c }
$\C W$ & \multicolumn{2}{|c|}{$\cdt{n}$} &  & \ $\csbt{n}{m}$ & & \multicolumn{2}{c}{$\csdt{n}{m}$} \\ \hline
& & & & & & & \\
Subspace & $\csdt{k+1}{1}$ &  $\cbt{k+1}$ & $\csbt{k}{m-r+1}$ & $\csbt{k+1}{r}$ & $\cct{k+1}$ &  $\csdt{k+1}{m-r+1}$ & $\csbt{k+1}{r}$ \\
& & & & & & & \\
Stabilizer & $A_{n-k-1}$ & $D_{n-k}$ &  $A_{n-k-1}$ &  $B_{n-k-1}$ &  $D_{{n-k-1},{r}}$ &  $A_{n-k-1}$ &  $D_{{n-k-1},{r}}$ \\
& & & & & & & \\
Enumeration & $2^{n-k-1} \binom{n}{n-k}$ &  $2 \binom{n}{n-k}$ &  $2^{n-k-1}\binom{n}{n-k}$ & $\binom{n-m}{n-k-r-1} \binom{m}{r}$ &
$\binom{n-m}{n-k-r-1} \binom{m}{r}$  & $2^{n-k-1}\binom{n}{n-k}$  &  $2 \binom{n-m}{n-k-r-1} \binom{m}{r}$ \\
& & & & & & & \\
\raisebox{.4in}{Configuration} &  \includegraphics{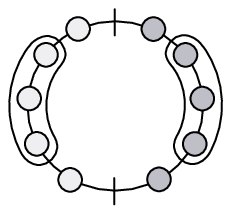} &  \includegraphics{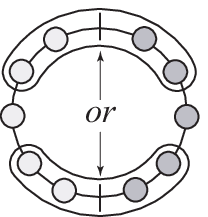} &
\includegraphics{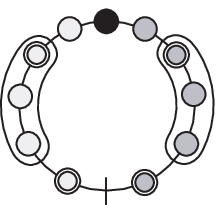} &  \includegraphics{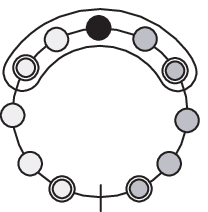} &  \includegraphics{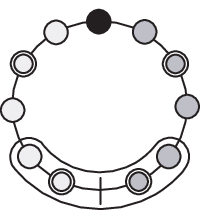} &
\includegraphics{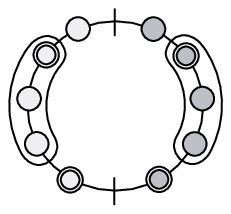} &  \includegraphics{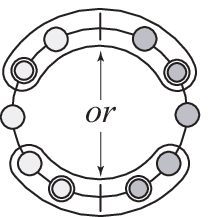} \\
\end{tabular}
\end{center}
\bigskip
\caption{Minimal building sets of dimension $k$ for Euclidean complexes.}
\label{t:euclid}
\end{table}
}
\end{landscape}

\clearpage

%%%%%%%%%%%%%%%%%%%%%%%%%%%%%%%%%%%%%%%%%%%%%%%%%%%%%%%%%%%%%%%%%%%%%%%%%%%%%%%%%%%%%
%
%                  REFERENCES
%
%%%%%%%%%%%%%%%%%%%%%%%%%%%%%%%%%%%%%%%%%%%%%%%%%%%%%%%%%%%%%%%%%%%%%%%%%%%%%%%%%%%%%
\bibliographystyle{amsplain}

\end{document}